\documentclass[12pt,leqno]{article}
\usepackage{amssymb}
\usepackage[mathscr]{eucal}
\usepackage{amsmath,amssymb,latexsym,theorem,bbm}
\usepackage{color}

\newcommand{\SC}{\scriptstyle}

\newcommand{\CC}{\mathsf{C}}
\newcommand{\DD}{\mathsf{D}}

\newcommand{\NN}{\mathbb{N}}

\newcommand{\RR}{\mathbb{R}}

\newcommand{\ZZ}{\mathbb{Z}}

\newcommand{\bA}{{\boldsymbol{A}}}

\newcommand{\bB}{{\boldsymbol{B}}}

\newcommand{\bC}{{\boldsymbol{C}}}

\newcommand{\bD}{{\boldsymbol{D}}}

\newcommand{\bI}{{\boldsymbol{I}}}

\newcommand{\bm}{{\boldsymbol{m}}}
\newcommand{\tbm}{\widetilde{\bm}}
\newcommand{\bM}{{\boldsymbol{M}}}
\newcommand{\bP}{{\boldsymbol{P}}}

\newcommand{\bq}{{\boldsymbol{q}}}

\newcommand{\bcQ}{{\boldsymbol{\cQ}}}
\newcommand{\bcR}{{\boldsymbol{\cR}}}
\newcommand{\bu}{{\boldsymbol{u}}}
\newcommand{\tbu}{\widetilde{\bu}}
\newcommand{\bv}{{\boldsymbol{v}}}
\newcommand{\bV}{{\boldsymbol{V}}}
\newcommand{\tbv}{\widetilde{\bv}}
\newcommand{\bx}{{\boldsymbol{x}}}
\newcommand{\bX}{{\boldsymbol{X}}}
\newcommand{\by}{{\boldsymbol{y}}}

\newcommand{\bU}{{\boldsymbol{U}}}

\newcommand{\bgamma}{{\boldsymbol{\gamma}}}
\newcommand{\bxi}{{\boldsymbol{\xi}}}

\newcommand{\balpha}{{\boldsymbol{\alpha}}}

\newcommand{\bvare}{{\boldsymbol{\vare}}}
\newcommand{\bPi}{{\boldsymbol{\Pi}}}

\newcommand{\bzero}{{\boldsymbol{0}}}

\newcommand{\cA}{{\mathcal A}}

\newcommand{\cD}{{\mathcal D}}

\newcommand{\cF}{{\mathcal F}}

\newcommand{\cM}{{\mathcal M}}
\newcommand{\bcM}{\boldsymbol{\cM}}

\newcommand{\cN}{{\mathcal N}}
\newcommand{\bcN}{\boldsymbol{\cN}}
\newcommand{\cP}{{\mathcal P}}
\newcommand{\cQ}{{\mathcal Q}}
\newcommand{\cR}{{\mathcal R}}

\newcommand{\cU}{{\mathcal U}}

\newcommand{\bcU}{\boldsymbol{\cU}}

\newcommand{\cX}{{\mathcal X}}
\newcommand{\bcX}{\boldsymbol{\cX}}

\newcommand{\cY}{{\mathcal Y}}
\newcommand{\cZ}{{\mathcal Z}}
\newcommand{\cW}{{\mathcal W}}
\newcommand{\bcW}{\boldsymbol{\cW}}

\newcommand{\bcY}{\boldsymbol{\cY}}

\newcommand{\dd}{\mathrm{d}}

\newcommand{\slu}{{\SC\mathrm{lu}}}

\newcommand{\INARp}{\textup{INAR($p$)}}

\newcommand{\EE}{\operatorname{\mathbb{E}}}
\newcommand{\PP}{\operatorname{\mathbb{P}}}
\newcommand{\OO}{\operatorname{O}}
\newcommand{\tr}{\operatorname{tr}}
\newcommand{\var}{\operatorname{Var}}

\newcommand{\varr}{\varrho}
\newcommand{\vare}{\varepsilon}

\renewcommand{\mid}{\,|\,}
\newcommand{\bmid}{\,\big|\,}

\renewcommand{\leq}{\leqslant}
\renewcommand{\geq}{\geqslant}

\newcommand{\stoch}{\stackrel{\PP}{\longrightarrow}}
\newcommand{\distr}{\stackrel{\cD}{\longrightarrow}}
\newcommand{\distre}{\stackrel{\cD}{=}}

\newcommand{\lu}{\stackrel{\slu}{\longrightarrow}}

\newcommand{\bbone}{\mathbbm{1}}
\newcommand{\ns}{{\lfloor ns\rfloor}}
\newcommand{\nt}{{\lfloor nt\rfloor}}
\newcommand{\nT}{{\lfloor nT\rfloor}}
\newcommand{\proofend}{\hfill\mbox{$\Box$}}

\numberwithin{equation}{section}

\theoremstyle{change} \theorembodyfont{\em}
\newtheorem{Lem}{Lemma.}[section]
\newtheorem{Thm}{Theorem.}[section]

\theorembodyfont{\rm}
\newtheorem{Rem}{Remark.}[section]

\begin{document}

\begin{center}
 {\bfseries\Large Asymptotic behavior of critical indecomposable} \\[2mm]
 {\bfseries\Large multi-type branching processes with immigration} \\[5mm]

 {\sc\large Tivadar Danka} and {\sc\large Gyula Pap}
\end{center}

\noindent
 Bolyai Institute, University of Szeged,
 Aradi v\'ertan\'uk tere 1, H--6720 Szeged, Hungary.

\noindent e--mails: theodore.danka@gmail.com (T. Danka),
                    papgy@math.u-szeged.hu (G. Pap).

\vskip-0.2cm

\footnote{The research of T. Danka and G. Pap was realized in the frames of
 T\'AMOP 4.2.4.\ A/2-11-1-2012-0001 'National Excellence Program --
 Elaborating and operating an inland student and researcher personal support
 system'.
The project was subsidized by the European Union and co-financed by the
 European Social Fund.}


\vspace*{-5mm}

\begin{abstract}
In this paper the asymptotic behavior of a critical multi-type branching
 process with immigration is described when the offspring mean matrix is
 irreducible, in other words, when the process is indecomposable. 
It is proved that sequences of appropriately scaled random step functions
 formed from periodic subsequences of a critical indecomposable multi-type
 branching process with immigration converge weakly towards a process
 supported by a ray determined by the Perron vector of the offspring mean
 matrix. 
The types can be partitioned into nonempty mutually disjoint subsets
 (according to communication of types) such that the coordinate processes
 belonging to the same subset are multiples of the same squared Bessel
 process, and the coordinate processes belonging to different subsets are
 independent.
\end{abstract}

\section{Introduction}
\label{section_intro}

Branching processes have a number of applications in biology, finance,
 economics, queueing theory etc., see e.g.\ Haccou, Jagers and Vatutin
 \cite{HJV}.
Many aspects of applications in epidemiology, genetics and cell kinetics were
 presented at the 2009 Badajoz Workshop on Branching Processes, see \cite{V}.

Let \ $(X_k)_{k\in\ZZ_+}$ \ be a single-type Galton--Watson branching process
 with immigration given by
 \begin{equation}\label{MBPI(1)}
  X_k = \sum_{j=1}^{X_{k-1}} \xi_{k, j} + \vare_k , \qquad k \in \NN , 
 \end{equation}
 and with initial value \ $X_0 = 0$.
\ Suppose that \ $\EE(\xi_{1,1}^2) < \infty$, \ $\EE(\vare_1^2) < \infty$, \ and
 \ $m_\xi := \EE(\xi_{1,1}) = 1$, \ i.e., the process is critial.
Wei and Winnicki \cite{WW} proved a functional limit theorem  
 \begin{equation}\label{conv_X}
  \cX^{(n)} \distr \cX \qquad \text{as \ $n \to \infty$,}
 \end{equation}
 where \ $\cX^{(n)}_t := n^{-1} X_{\nt}$ \ for \ $t \in \RR_+$, \ $n \in \NN$,
 \ and \ $(\cX_t)_{t\in\RR_+}$ \ is the pathwise unique strong solution of the
 stochastic differential equation (SDE)
 \[
   \dd \cX_t = m_\vare \, \dd t + \sqrt{ V_\xi \cX_t^+ } \, \dd \cW_t ,
   \qquad t \in \RR_+ ,
 \]
 with initial value \ $\cX_0 = 0$, \ where \ $m_\vare := \EE(\vare_1)$,
 \ $V_\xi := \var(\xi_{1,1})$, \ $(\cW_t)_{t\in\RR_+}$ \ is a standard Wiener
 process, and \ $x^+$ \ denotes the positive part of \ $x \in \RR$.

A multi-type branching process \ $(\bX_k)_{k\in\ZZ_+}$ \ is referred to
 respectively as subcritical, critical or supercritical if
 \ $\varr(\bm_\bxi) < 1$, \ $\varr(\bm_\bxi) = 1$ \ or \ $\varr(\bm_\bxi) > 1$,
 \ where \ $\varr(\bm_\bxi)$ \ denotes the spectral radius of the offspring
 mean matrix \ $\bm_\bxi$ \ (see, e.g., Athreya and Ney \cite{AN} or
 Quine \cite{Q}).
Joffe and M\'etivier \cite[Theorem 4.3.1]{JM} studied a critical multi-type
 branching process without immigration when the offspring mean matrix is
 primitive, in other words, when the process is positively regular.
They determined the limiting behavior of the martingale part
 \ $(\bcM^{(n)})_{n \in \NN}$ \ given by
 \ $\bcM^{(n)}_t := n^{-1} \sum_{k=1}^{\nt} \bM_k$ \ with
 \ $\bM_k := \bX_k - \EE(\bX_k \mid \bX_0, \dots, \bX_{k-1})$ \ (which is a
 special case of Theorem \ref{Conv_M}). 

The result \eqref{conv_X} has been generalized by Isp\'any and Pap \cite{IP}
 for a critical $p$-type branching process \ $(\bX_k)_{k\in\ZZ_+}$ \ with
 immigration when the offspring mean matrix is primitive, in other words, when
 the process is positively regular.
They proved that
 \[
   \bcX^{(n)} \distr \cX \bu \qquad \text{as \ $n \to \infty$,}
 \]
 where \ $\bcX^{(n)}_t := n^{-1} \bX_{\nt}$ \ for \ $t \in \RR_+$, \ $n \in \NN$,
 \ the process \ $(\cX_t)_{t\in\RR_+}$ \ is the unique strong solution of the SDE
 \[
   \dd \cX_t
   = \bv^\top \! \bm_\vare \, \dd t
     + \sqrt{\bv^\top (\bu \odot \bV_\bxi) \bv \cX_t^+} \, \dd \cW_t ,
   \qquad t \in \RR_+ ,
 \]
 with initial value \ $\cX_0 = 0$, \ where \ $(\cW_t)_{t\in\RR_+}$ \ is a
 standard Wiener process, \ $\bu$ \ and \ $\bv$ \ denotes the right and left
 Perron vectors of \ $\bm_\bxi$, \ $\bm_\vare$ \ denotes the immigration mean
 vector, and \ $\bu \odot \bV_\bxi := \sum_{i=1}^p u_i \bV_{\bxi_i}$, \ where
 \ $\bu = (u_i)_{i\in\{1,\ldots,p\}}$ \ and \ $\bV_{\bxi_i}$ \ denotes the offspring
 variance matrix of type \ $i \in \{1, \ldots, p\}$.

The aim of the present paper is to obtain a generalization of this result for a
 critical multi-type branching process with immigration when the offspring mean
 matrix is irreducible, in other words, when the process is indecomposable.
We succeeded to determine the joint asymptotic behavior of sequences
 \ $\bigl((nr)^{-1} \bX_{r\nt+i-1}\bigr)_{t\in\RR_+}$, \ $n \in \NN$,
 \ $i \in \{1, \ldots, r\}$, \ of random step functions as \ $n \to \infty$,
 \ where \ $r$ \ denotes the index of cyclicity (in other words, the index of
 imprimivity) of the offspring mean matrix (see Theorem \ref{main}).
It turns out that the limiting diffusion process has the form
 \ $(\bm_\bxi^{r-i+1} \bcY_t)_{t\in\RR_+}$, \ $i \in \{1, \ldots, r\}$, \ where the
 distribution of the process \ $(\bm_\bxi^j \bcY_t)_{t\in\RR_+}$ \ is the same for
 all \ $j \in \ZZ_+$.
\ Moreover, the process \ $(\bcY_t)_{t\in\RR_+}$ \ is 1-dimensional in the sense
 that for all \ $t \in \RR_+$, \ the distribution of \ $\bcY_t$ \ is
 concentrated on the ray \ $\RR_+ \cdot \bu$, \ where \ $\bu$ \ is the Perron
 eigenvector of the offspring mean matrix \ $\bm_\bxi$.
\ In fact, the types can be partitioned into nonempty mutually disjoint subsets
 (according to communication of types) such that partitioning the coordinates
 of the limit process \ $\bcY_t = (\bcY_{t,1}, \ldots, \bcY_{t,r})$, \ and of the
 Perron eigenvector \ $\bu = (\bu_1, \ldots, \bu_r)$, \ we have
 \ $\bcY_{t,i} = \cZ_{t,i} \, \bu_i$, \ $t \in \RR_+$,
 \ $i \in \{1, \ldots, r\}$, \ where \ $(\cZ_{t,i})_{t\in\RR_+}$,
 \ $i \in \{1, \ldots, r\}$, \ are independent squared Bessel processes.
It is interesting to note that Kesten and Stigum \cite{KS2} considered a
 supercritical indecomposable multi-type branching processes without
 immigration, and they proved that there exits a random variable \ $w$ \ such
 that \ $\varrho(\bm_\xi)^{-(r n + i - 1)} \bX_{r n + i - 1} \to w \bu_i$ \ almost
 surely as \ $n \to \infty$ \ for each \ $i \in \{1, \ldots, r\}$.

\section{Multi-type branching processes with immigration}

Let \ $\ZZ_+$, \ $\NN$, \ $\RR$, \ $\RR_+$  \ and \ $\RR_{++}$ \ denote the set
 of non-negative integers, positive integers, real numbers, non-negative real
 numbers and positive real numbers, respectively.
Every random variable will be defined on a fixed probability space
 \ $(\Omega,\cA,\PP)$.

We will investigate a $p$-type branching process \ $(\bX_k)_{ k\in\ZZ_+ }$ \ with
 immigration.
For simplicity, we suppose that the initial value is \ $\bX_0 = \bzero$.
\ For each \ $k \in \ZZ_+$ \ and \ $i \in \{ 1, \dots, p \}$, \ the number of
 individuals of type \ $i$ \ in the \ $k^\mathrm{th}$ \ generation is denoted by
 \ $X_{k,i}$.
\ The non-negative integer-valued random variable \ $\xi_{k,j,i,\ell}$ \ denotes
 the number of type \ $\ell$ \ offsprings produced by the \ $j^\mathrm{th}$
 \ individual who is of type \ $i$ \ belonging to the \ $(k-1)^\mathrm{th}$
 \ generation.
The number of type \ $i$ \ immigrants in the \ $k^\mathrm{th}$ \ generation will
 be denoted by \ $\vare_{k,i}$.
\ Consider the random vectors
 \[
   \bX_k := \begin{bmatrix}
                  X_{k,1} \\
                  \vdots \\
                  X_{k,p}
                 \end{bmatrix} , \qquad
   \bxi_{k,j,i} := \begin{bmatrix}
                  \xi_{k,j,i,1} \\
                  \vdots \\
                  \xi_{k,j,i,p}
                 \end{bmatrix} , \qquad
   \bvare_k := \begin{bmatrix}
                \vare_{k,1} \\
                \vdots \\
                \vare_{k,p}
               \end{bmatrix} .
 \]
Then we have
 \begin{equation}\label{MBPI(d)}
  \bX_k = \sum_{i=1}^p \sum_{j=1}^{X_{k-1,i}} \bxi_{k,j,i} + \bvare_k , \qquad
  k \in \NN .
 \end{equation}
Here
 \ $\big\{\bxi_{k,j,i}, \, \bvare_k
          : k, j \in \NN, \, i \in \{ 1, \dots, p \} \big\}$
 \ are supposed to be independent. 
Moreover, \ $\big\{ \bxi_{k,j,i} : k, j \in \NN \big\}$ \ for each
 \ $i \in \{ 1, \dots, p \}$, \ and \ $\{ \bvare_k : k \in \NN \}$ \ are
 supposed to consist of identically distributed vectors.

Suppose \ $\EE(\|\bxi_{1,1,i}\|^2) < \infty$ \ for all
 \ $i \in \{ 1, \dots, p \}$ \ and \ $\EE (\|\bvare_1\|^2) < \infty$.
\ Introduce the notations
 \begin{gather*}
  \bm_{\bxi_i} := \EE(\bxi_{1,1,i}) \in \RR^p_+ , \qquad
  \bm_{\bxi} := \begin{bmatrix}
                \bm_{\bxi_1} & \cdots & \bm_{\bxi_d}
               \end{bmatrix} \in \RR^{p \times p}_+ , \qquad
  \bm_{\bvare} := \EE(\bvare_1) \in \RR^p_+ , \\
  \bV_{\bxi_i} := \var(\bxi_{1,1,i}) \in \RR^{p \times p} , \qquad
  \bV_{\bvare} := \var(\bvare_1) \in \RR^{p \times p} .
 \end{gather*}
Note that many authors define the offspring mean matrix as \ $\bm^\top_\bxi$. 
\ For \ $k \in \ZZ_+$, \ let \ $\cF_k := \sigma(\bX_0 ,\bX_1 , \dots, \bX_k)$.
\ By \eqref{MBPI(d)}, 
 \begin{equation}\label{mart}
  \EE(\bX_k \mid \cF_{k-1}) 
  = \sum_{i=1}^p X_{k-1,i} \bm_{\bxi_i} + \bm_{\bvare}
  = \bm_{\bxi} \bX_{k-1} + \bm_{\bvare} .
 \end{equation}
Consequently,
 \begin{equation}\label{recEX}
  \EE(\bX_k) 
  = \bm_{\bxi} \EE(\bX_{k-1}) + \bm_{\bvare} , \qquad
  k \in \NN ,
 \end{equation}
 which implies
 \begin{equation}\label{EXk}
  \EE(\bX_k) 
  = \sum_{j=0}^{k-1} \bm_{\bxi}^j \bm_{\bvare} , \qquad
  k \in \NN .
 \end{equation}
Hence the offspring mean matrix \ $\bm_{\bxi}$ \ plays a crucial role in the
 asymptotic behavior of the sequence \ $(\bX_k)_{k\in\ZZ_+}$.

In what follows we recall some known facts about irreducible nonnegative
 matrices.
The matrix \ $\bm_{\bxi}$ \ is reducible if there exist a permutation matrix
 \ $\bP \in \RR^{p \times p}$ \ and an integer \ $q$ \ with \ $1 \leq q \leq p-1$
 \ such that
 \[
   \bP^\top \bm_{\bxi} \bP
   = \begin{bmatrix} \bB & \bC \\ \bzero & \bD \end{bmatrix} ,
 \]
 where \ $\bB \in \RR^{q \times q}$, \ $\bD \in \RR^{(p-q) \times (p-q)}$, \ 
 \ $\bC \in \RR^{q \times (p-q)}$ \ and \ $\bzero \in \RR^{(p-q) \times q}$ \ is a
 null matrix.
The matrix \ $\bm_{\bxi}$ \ is irreducible if it is not reducible; see, e.g.,
 Horn and Johnson \cite[Definition 6.2.21, Definition 6.2.22]{HJ}.
The matrix \ $\bm_{\bxi}$ \ is irreducible if and only if for all
 \ $i, j \in \{1, \ldots, p\}$ \ there exists \ $n_{i,j} \in \NN$ \ such that
 the matrix entry \ $(\bm_\bxi^{n_{i,j}})_{i,j}$ \ is positive.
If the matrix \ $\bm_{\bxi}$ \ is irreducible then, by the Frobenius--Perron
 theorem the following assertions hold
 (see, e.g., Bapat and Raghavan \cite[Theorem 1.8.3]{BR},
  Berman and Plemmons \cite[Theorem 2.20]{BR},
  Brualdi and Cvetkovi\'c \cite[Theorem 8.2.4]{BC}, Minc \cite[Theorem 3.1]{M}
  or Kesten and Stigum \cite{KS2}) :
 \begin{itemize}
  \item
   $\varrho(\bm_{\bxi}) \in \RR_{++}$, \ $\varrho(\bm_{\bxi})$ \ is an eigenvalue
    of \ $\bm_{\bxi}$, \ and the algebraic and geometric multiplicities of
    \ $\varrho(\bm_{\bxi})$ \ equal 1.
  \item
   Corresponding to the eigenvalue \ $\varrho(\bm_{\bxi})$ \ there exists a
    unique (right) eigenvector \ $\bu \in \RR^p_{++}$, \ called the Perron
    vector of \ $\bm_{\bxi}$, \ such that the sum of the coordinates of
    \ $\bu$ \ is 1, and there exists a unique left eigenvector
    \ $\bv \in \RR^p_{++}$, \ such that \ $\bu^\top \bv = 1$.
  \item
   If \ $\bm_{\bxi}$ \ has exactly \ $r$ \ eigenvalues of maximum modulus
    \ $\varrho(\bm_{\bxi})$ \ then the coordinates \ $\{1, \ldots, p\}$ \ can be
    partitioned into \ $r$ \ nonempty mutually disjoint subsets
    \ $D_1, \ldots, D_r$ \ such that \ $m_{i,j} = 0$ \ unless there exists 
    \ $\ell \in \{1, \ldots, r\}$ \ with \ $i \in D_{\ell-1}$ \ and
    \ $j \in D_\ell$, \ where subscripts are considered modulo \ $r$.
   \ This partitioning is unique up to cyclic permutation of the subsets.
   We may assume that the types are enumerated according to these subsets, and
    hence
    \begin{equation}\label{partition_m}
      \bm_{\bxi}
      = \begin{bmatrix}
         \bzero   & \bm_{1,2} & \bzero   & \cdots & \bzero & \bzero \\
         \bzero   & \bzero   & \bm_{2,3} & \cdots & \bzero & \bzero \\
         \bzero   & \bzero   & \bzero   & \cdots & \bzero & \bzero \\
         \vdots   & \vdots   & \vdots   & \ddots & \vdots & \vdots \\
         \bzero   & \bzero   & \bzero   & \cdots & \bzero & \bm_{r-1,r} \\
         \bm_{r,1} & \bzero   & \bzero   & \cdots & \bzero & \bzero
        \end{bmatrix} ,
    \end{equation}
    where the $r$ main diagonal zero blocks are square,
    \ $\bm_{1,2} \in \RR^{|D_1| \times |D_2|}$,
    \ $\bm_{2,3} \in \RR^{|D_2| \times |D_3|}$,
    \ \dots, \ $\bm_{r,1} \in \RR^{|D_r| \times |D_1|}$ \ where \ $|H|$ \ denotes
    the number of elements of a set \ $H$, \ and \ $\bm_{1,2} \ne \bzero$,
    \ $\bm_{2,3} \ne \bzero$, \ \dots, \ $\bm_{r,1} \ne \bzero$.
   \ Then for each \ $k \in \{1, \ldots, r-1\}$, \ we have
    \[
      \bm_{\bxi}^k
      = \begin{bmatrix}
         \bzero & \cdots & \bzero & \tbm_{1,k+1} & \bzero   & \cdots & \bzero \\
         \bzero & \cdots & \bzero & \bzero & \tbm_{2,k+2} & \cdots & \bzero \\
         \vdots & & \vdots & \vdots & \vdots & \ddots & \vdots \\
         \bzero & \cdots & \bzero & \bzero & \bzero & \cdots & \tbm_{r-k,r} \\
         \tbm_{r-k+1,r+1} & \cdots & \bzero & \bzero & \bzero & \cdots & \bzero \\
         \vdots & \ddots & \vdots & \vdots & \vdots & & \vdots \\
         \bzero & \cdots & \tbm_{r,r+k} & \bzero & \bzero & \cdots & \bzero
        \end{bmatrix}
    \]
    with
    \[
      \tbm_{i,j}
      := \bm_{i,i+1} \bm_{i+1,i+2} \cdots \bm_{j-1,j} \in \RR^{|D_i| \times |D_j|}
    \]
    for \ $i, j \in \NN$ \ with \ $i < j$, \ where subscripts on the right
    hand side are considered modulo \ $r$.
   \ We will also use the convention that
    \ $\tbm_{i,i} \in \RR^{|D_i| \times |D_i|}$ \ is the identity matrix.
   Moreover,
    \[
      \bm_{\bxi}^r
      = \begin{bmatrix}
         \tbm_{1,r+1}  & \bzero & \cdots & \bzero \\
         \bzero & \tbm_{2,r+2}  & \cdots & \bzero \\
         \vdots & \vdots & \ddots & \vdots \\
         \bzero & \bzero & \cdots & \tbm_{r,2r}
        \end{bmatrix}
      =: \tbm_{1,r+1} \oplus \tbm_{2,r+2} \oplus \cdots \oplus \tbm_{r,2r}
    \]
   The matrices \ $\tbm_{1,r+1} \in \RR^{|D_1| \times |D_1|}$,
    \ $\tbm_{2,r+2} \in \RR^{|D_2| \times |D_2|}$, \ldots,
    $\tbm_{r,2r} \in \RR^{|D_r| \times |D_r|}$ \ are primitive (that is,
    irreducible and there exists \ $n_i \in \NN$ \ such that
    \ $\tbm_{i,i+r}^{n_i} \in \RR_{++}^{|D_i| \times |D_i|}$) with
    \ $\varrho(\tbm_{i,i+r}) = [\varrho(\bm_{\bxi})]^r$,
    \ $i \in \{1, \ldots, r\}$.
   \ (See, e.g., Minc \cite[Theorem 4.3]{M}.)
  \item
   If
    \[
      \bu = \begin{bmatrix}
             \bu_1 \\
             \vdots \\
             \bu_r \\
            \end{bmatrix} , \qquad
      \bv = \begin{bmatrix}
             \bv_1 \\
             \vdots \\
             \bv_r \\
            \end{bmatrix}
    \]
    denotes the partitioning of \ $\bu$ \ and \ $\bv$ \ with respect to
    the partitioning \ $D_1, \ldots, D_r$ \ of the coordinates
    \ $\{1, \ldots, p\}$ \ then, for each \ $i \in \{1, \ldots, r\}$,
    \ the vector \ $\tbu_i := r \bu_i$ \ is the Perron eigenvector of
    \ $\tbm_{i,i+r}$, \ and \ $\tbv_i := \bv_i$ \ is the corresponding left
    eigenvector.
   Further,
    \[
      [\varrho(\tbm_{i,i+r})]^{-n} \, \tbm_{i,i+r}^n
      \to \tbu_i \tbv_i^\top 
      = r \bu_i \bv_i^\top 
      =: \bPi_i \in \RR^{|D_i| \times |D_i|}_{++} \qquad
      \text{as \ $n \to \infty$.} 
    \]
   Consequently,
    \[
      [\varrho(\bm_{\bxi})]^{-nr} \, \bm_{\bxi}^{nr}
      \to \bPi \in \RR^{p \times p}_+ \qquad
      \text{as \ $n \to \infty$,} 
    \]
    where 
    \begin{equation}\label{bPi}
      \bPi
      := \bPi_1 \oplus \bPi_2 \oplus \cdots \oplus \bPi_r
       = \begin{bmatrix}
          \bPi_1  & \bzero  & \cdots & \bzero \\
          \bzero  & \bPi_2  & \cdots & \bzero \\
          \vdots  & \vdots  & \ddots & \vdots \\
          \bzero  & \bzero  & \cdots & \bPi_r \\
         \end{bmatrix} .
    \end{equation}
   The vectors \ $\bu$ \ and \ $\bv$ \ are right and left eigenvectors of
    \ $\bPi$ \ corresponding to the eigenvalue \ $[\varrho(\bm_{\bxi})]^r$.
  \item
   Moreover, there exist \ $c , \kappa \in \RR_{++}$ \ with \ $\kappa < 1$
    \ such that 
    \begin{equation}\label{rate}
     \| [\varrho(\bm_{\bxi})]^{-nr} \, \bm_{\bxi}^{nr} - \bPi \| \leq c \kappa^n
     \qquad \text{for all \ $n \in \NN$,}
    \end{equation}
    where \ $\|\bA\|$ \ denotes the operator norm of a matrix
    \ $\bA \in \RR^{p \times p}$ \ defined by
    \ $\|\bA\| := \sup_{\|\bx\| = 1} \| \bA \bx \|$.
 \end{itemize}
The number \ $r$ \ of eigenvalues of maximum modulus \ $\varrho(\bm_{\bxi})$
 \ is called the index of cyclicity (in other words, the index of imprimivity)
 of the matrix \ $\bm_{\bxi}$.

If \ $\bm_\bxi$ \ has the form \eqref{partition_m}, then the offsprings have
 the property \ $\xi_{1,1,i,j} = 0$ \ $\PP$-almost surely unless there exists
 \ $\ell \in \{1, \ldots, r\}$ \ with \ $i \in D_{\ell-1}$ \ and
 \ $j \in D_\ell$, \ where subscripts considered modulo \ $r$.
\ Consequently, the offspring variance matrices \ $\bV_{\bxi_j}$,
 \ $j \in \{1, \ldots, p\}$, \ have the form
 \begin{equation}\label{Cov_V}
  \bV_{\bxi_j}
  = \begin{cases}
     \bzero \oplus \bzero \oplus \cdots \oplus \bzero \oplus \bV_{1,j}
      & \text{if \ $j \in D_1$,} \\
     \bV_{2,j} \oplus \bzero \oplus \cdots \oplus \bzero \oplus \bzero
      & \text{if \ $j \in D_2$,} \\
     \vdots \\
     \bzero \oplus \bzero \oplus \cdots \oplus \bV_{r,j} \oplus \bzero
      & \text{if \ $j \in D_r$,}
    \end{cases}
 \end{equation}
 where \ $\bV_{\ell,j} \in \RR^{|D_{\ell-1}| \times |D_{\ell-1}|}$ \ denotes the
 variance matrix of the random vector \ $(\xi_{1,1,j,i})_{i \in D_{\ell-1}}$ \ for
 \ $\ell \in \{1, \ldots, r\}$, \ $j \in D_\ell$, \ where subscripts considered
 modulo \ $r$.
\ For a vector
 \ $\balpha_\ell = (\alpha_{\ell,j})_{j \in D_\ell} \in \RR^{|D_\ell|}_+$ \ with
 \ $\ell \in \{1, \ldots, r\}$, \ we will use notation
 \ $\balpha_\ell \odot \bV_\ell := \sum_{j \in D_\ell} \alpha_{\ell,j} \bV_{\ell,j}
    \in \RR^{|D_{\ell-1}| \times |D_{\ell-1}|}$,
 \ which is a positive semi-definite matrix, a mixture of the variance
 matrices \ $\bV_{\ell,j}$, \ $j \in D_\ell$.
\ For a vector \ $\balpha = (\alpha_i)_{i=1,\dots,p} \in \RR^p_+$, \ we will also
 use the notation
 \ $\balpha \odot \bV_\bxi
    := \sum_{i=1}^p \alpha_i \bV_{\bxi_i} \in \RR^{p \times p}$,
 \ which is a positive semi-definite matrix, a mixture of the variance
 matrices \ $\bV_{\xi_1}, \ldots, \bV_{\xi_p}$.

A multi-type branching process with immigration is called indecomposable if its
 offspring mean matrix \ $\bm_\bxi$ \ is irreducible.

\section{Convergence results}
\label{Convergence_result}

A function \ $f : \RR_+ \to \RR^p$ \ is called \emph{c\`adl\`ag} if it is right
 continuous with left limits.
\ Let \ $\DD(\RR_+, \RR^p)$ \ and \ $\CC(\RR_+, \RR^p)$ \ denote the space of
 all $\RR^p$-valued c\`adl\`ag and continuous functions on \ $\RR_+$,
 \ respectively.
Let \ $\cD_\infty(\RR_+, \RR^p)$ \ denote the Borel $\sigma$-algebra in
 \ $\DD(\RR_+, \RR^p)$ \ for the metric defined in Jacod and Shiryaev
 \cite[Chapter VI, (1.26)]{JSh} (with this metric \ $\DD(\RR_+, \RR^p)$ \ is a
 complete and separable metric space).
For $\RR^p$-valued stochastic processes \ $(\bcY_t)_{t\in\RR_+}$ \ and
 \ $(\bcY^{(n)}_t)_{t\in\RR_+}$, \ $n \in \NN$, \ with c\`adl\`ag paths we write
 \ $\bcY^{(n)} \distr \bcY$ \ if the distribution of \ $\bcY^{(n)}$ \ on the
 space \ $(\DD(\RR_+, \RR^p), \cD_\infty(\RR_+, \RR^p))$ \ converges weakly to
 the distribution of \ $\bcY$ \ on the space
 \ $(\DD(\RR_+, \RR^p), \cD_\infty(\RR_+, \RR^p))$ \ as \ $n \to \infty$.

For each \ $n \in \NN$, \ consider the \ $rp$-dimensional random step process
 \[
   \bcX^{(n)}_t := (nr)^{-1} \begin{bmatrix}
                            \bX_{r\nt} \\
                            \bX_{r\nt-1} \\
                            \vdots \\
                            \bX_{r\nt-r+1}
                           \end{bmatrix} ,
   \qquad t \in \RR_+ ,
 \]
 where \ $\lfloor x \rfloor$ \ denotes the integer part of \ $x \in \RR$. 
The positive part of \ $x \in \RR$ \ will be denoted by \ $x^+$.

\begin{Thm}\label{main}
Let \ $(\bX_k)_{k\in\ZZ_+}$ \ be a critical indecomposable $p$-type branching
 process with immigration.
Suppose \ $\bX_0 = \bzero$, \ $\EE(\|\bxi_{1,1,i}\|^4) < \infty$ \ for all
 \ $i \in \{ 1, \dots, p \}$ \ and \ $\EE (\|\bvare_1\|^4) < \infty$.
\ Suppose that the offspring mean matrix \ $\bm_\bxi$ \ has the form
 \eqref{partition_m}.
Then
 \begin{gather}\label{Conv_X}
  \bcX^{(n)} \distr \bcX \qquad \text{as \ $n \to \infty$,}
 \end{gather}
 where
 \[
   \bcX_t = \begin{bmatrix}
             \bm_\bxi^r \bcY_t \\
             \bm_\bxi^{r-1} \bcY_t \\
             \vdots \\
             \bm_\bxi \bcY_t
            \end{bmatrix}
          = \begin{bmatrix}
             \bcY_t \\
             \bm_\bxi^{r-1} \bcY_t \\
             \vdots \\
             \bm_\bxi \bcY_t
            \end{bmatrix} , \qquad t \in \RR_+ ,
 \]
 with
 \[
   \bcY_t = \begin{bmatrix}
             \bcY_{t,1} \\
             \bcY_{t,2} \\
             \vdots \\
             \bcY_{t,r}
            \end{bmatrix} , \qquad t \in \RR_+ ,
 \]
 where, for \ $i \in \{1, \ldots, r\}$, \ the $|D_i|$-dimensional process
 \ $(\bcY_{t,i})_{t\in\RR_+}$ \ is given by
 \[
   \bcY_{t,i} := \cZ_{t,i} \bu_i , \qquad t \in\RR_+ ,
 \]
 where \ $(\cZ_{t,i})_{t\in\RR_+}$ \ is the unique strong solution of the SDE 
 \begin{equation}\label{SDE_Z}
  \dd \cZ_{t,i}
  = \bv_i^\top \bm_{\bxi,\bvare,i} \, \dd t
    + \sqrt{ \bv_i^\top \bV_{\bxi,\bvare,i} \bv_i \cZ_{t,i}^+ } \, \dd \cW_{t,i} ,
  \qquad t \in \RR_+ ,
 \end{equation}
 with initial value \ $\cZ_{0,i} = 0$, \ where \ $(\cW_{t,i})_{t\in\RR_+}$,
 \ $i \in \{1, \ldots, r\}$, \ are independent standard Wiener processes and
 \[
   \bm_{\bxi,\bvare,i}
   := \frac{1}{r} \sum_{\ell=i}^{i+r-1} \tbm_{i,\ell} \bm_{\bvare,\ell} , \qquad
   \bV_{\bxi,\bvare,i}
   := \frac{1}{r}
      \sum_{\ell=i}^{i+r-1} 
       \tbm_{i,\ell}
       \left[ \left( \tbm_{\ell+1,i+r} \bu_i \right) \odot \bV_{\ell+1} \right]
       \tbm_{i,\ell}^\top ,
 \]
 where
 \[
   \bm_\bvare = \begin{bmatrix}
                \bm_{\bvare,1} \\
                \bm_{\bvare,2} \\
                \vdots \\
                \bm_{\bvare,r}
               \end{bmatrix}
 \]
 denotes the partitioning of \ $\bm_\bvare$ \ with respect to the partitioning
 \ $D_1, \ldots, D_r$ \ of the types \ $\{1, \ldots, p\}$.
\ Moreover, the $p$-dimensional coordinate blocks of the $rp$-dimensional
 process \ $(\bcX_t)_{t\in\RR_+}$ \ have the same distribution, i.e.,
 \ $(\bm_\bxi^i \bcY_t)_{t\in\RR_+} \distre (\bcY_t)_{t\in\RR_+}$ \ for
 all \ $i \in \{1, \ldots, r-1\}$, \ and they are periodic, i.e.,
 \ $(\bm_\bxi^r \bcY_t)_{t\in\RR_+} = (\bcY_t)_{t\in\RR_+}$ \ with probability one.
\end{Thm}

\begin{Rem}
The higher moment assumptions \ $\EE(\|\bxi_{1,1,i}\|^4) < \infty$,
 \ $i \in \{1, \dots, p\}$ \ and \ $\EE (\|\bvare_1\|^4) < \infty$ \ are
 needed only for checking the conditional Lindeberg condition, namely,
 condition (ii) of Theorem \ref{Conv2DiffThm} for proving convergence of the
 martingale part described in Theorem \ref{Conv_M}. 
We suspect that one can check the conditional Lindeberg condition under the
 weaker moment assumptions \ $\EE(\|\bxi_{1,1,i}\|^2) < \infty$,
 \ $i \in \{1, \dots, p\}$ \ and \ $\EE (\|\bvare_1\|^2) < \infty$ \ by the
 method of Isp\'any and Pap \cite{IP}, see also this method in Barczy et
 al.~\cite{BarIspPap0}.
\proofend
\end{Rem}

\begin{Rem}\label{Rem_SDE_Z}
The SDE \eqref{SDE_Z} has a unique strong solution
 \ $(\cZ_{t,\ell}^{(z_0)})_{t\in\RR_+}$ \ for all initial values
 \ $\cZ_{0,\ell}^{(z_0)} = z_0 \in \RR$.
\ Indeed, the coefficient functions satisfy conditions of part (ii) of
 Theorem 3.5 in Chapter IX in Revuz and Yor \cite{RY} or the conditions of
 Proposition 5.2.13 in Karatzas and Shreve \cite{KarShr}.
Further, by the comparison theorem
 (see, e.g., Revuz and Yor \cite[Theorem 3.7, Chapter IX]{RY}), if the
 initial value \ $\cZ_{0,\ell}^{(z_0)} = z_0$ \ is nonnegative, then
 \ $\cZ_{t,\ell}^{(z_0)}$ \ is nonnegative for all \ $t \in \RR_+$ \ with
 probability one.
Hence \ $\cZ_{t,\ell}^+$ \ may be replaced by \ $\cZ_{t,\ell}$ \ under the square
 root in \eqref{SDE_Z}.
\proofend
\end{Rem}

\begin{Rem}
Note that Theorem \ref{main} implies the convergence result of Isp\'any and
 Pap \cite{IP} for a critical $p$-type branching process \ $(\bX_k)_{k\in\ZZ_+}$
 \ with immigration when the offspring mean matrix is primitive.
Indeed, in this case the index of cyclicity is \ $r = 1$, \ and
 \ $\bm_{\bxi,\bvare,i} = \bm_{\bvare,\ell}$,
 \ $\bV_{\bxi,\bvare,i} = \bu_i \odot \bV_\bxi$.
\proofend
\end{Rem}

In order to prove \eqref{Conv_X}, introduce the $rp$-dimensional random
 vectors
 \begin{equation}\label{Mk}
  \bM_k
  := \begin{bmatrix}
      \bM_{k,1} \\
      \bM_{k,2} \\
      \vdots \\
      \bM_{k,r} \\
     \end{bmatrix}
  := \begin{bmatrix}
      \bX_{rk} - \EE(\bX_{rk} \mid \cF_{rk-1}) \\
      \bX_{rk-1} - \EE(\bX_{rk-1} \mid \cF_{rk-2}) \\
      \vdots \\
      \bX_{rk-r+1} - \EE(\bX_{rk-r+1} \mid \cF_{rk-r}) \\
     \end{bmatrix}
   = \begin{bmatrix}
      \bX_{rk} - \bm_{\bxi} \bX_{rk-1} - \bm_{\bvare} \\
      \bX_{rk-1} - \bm_{\bxi} \bX_{rk-2} - \bm_{\bvare} \\
      \vdots \\
      \bX_{rk-r+1} - \bm_{\bxi} \bX_{rk-r} - \bm_{\bvare}
     \end{bmatrix}
 \end{equation}
 for \ $k \in \NN$, \ forming a sequence of martingale differences with respect
 to the filtration \ $(\cF_{rk})_{k\in\ZZ_+}$.
\ Consider the $rp$-dimensional random step processes
 \[
   \bcM^{(n)}_t := \begin{bmatrix}
                   \bcM^{(n)}_{t,1} \\
                   \vdots \\
                   \bcM^{(n)}_{t,r} \\
                  \end{bmatrix}
               := (nr)^{-1} \sum_{k=1}^{\nt} \bM_k ,
   \qquad t \in \RR_+ , \qquad n \in \NN .
 \]
The following convergence result is an important step in the proof of Theorem
 \ref{main}.

\begin{Thm}\label{Conv_M}
Under the assumptions of Theorem \ref{main}, we have
 \[
   \bcM^{(n)} \distr \bcM \qquad \text{as \ $n \to \infty$,}
 \]
 where
 \[
   \bcM_t := \begin{bmatrix}
              \bcM_{t,1} \\
              \vdots \\
              \bcM_{t,r} \\
             \end{bmatrix}
   \qquad t \in \RR_+ ,
 \]
 is the unique strong solution of the SDE
 \begin{equation}\label{SDE_M}
  \dd \bcM_{t,i}
  = \frac{1}{r}
    \sqrt{ \left[ \bm_\bxi^{r-i} \bPi 
                  \sum_{j=1}^r
                   \bm_\bxi^{j-1} (r \bcM_{t,j} + t \bm_\bvare) \right]^+
            \odot \bV_\bxi } \,
    \dd \bcW_{t,i} , \qquad i \in \{1, \ldots, r\} ,
 \end{equation}
 with initial value \ $\bcM_0 = \bzero$, \ where \ $(\bcW_{t,i})_{t\in\RR_+}$,
 \ $i \in \{ 1, \ldots, r \}$, \ are independent standard $p$-dimensional
 Wiener processes, and for a positive semi-definite matrix
 \ $\bA \in \RR^{p \times p}$, \ $\sqrt{\bA}$ \ denotes its unique symmetric
 positive semi-definite square root.
\end{Thm}

In order to handle the SDE \eqref{SDE_M}, consider the $p$-dimensional process
 \begin{equation}\label{N}
  \bcN_t := \bcM_{t,1} + \bm_\bxi \bcM_{t,2} + \cdots + \bm_\bxi^{r-1} \bcM_{t,r} ,
  \qquad t \in \RR_+ .
 \end{equation}

\begin{Thm}\label{PROC_N}
Under the assumptions of Theorem \ref{main}, the process \ $(\bcN_t)_{t\in\RR_+}$
 \ is the unique strong solution of the SDE
 \begin{equation}\label{SDE_N}
  \dd \bcN_t
  = \frac{1}{r}
    \sum_{j=1}^r
     \bm_\bxi^{j-1}
     \sqrt{ \left[ \bm_\bxi^{r-j} \bPi
                   \left(r \bcN_t
                         + t \sum_{\ell=1}^r \bm_\bxi^{\ell-1} \bm_\bvare\right)
                   \right]^+
            \odot \bV_\bxi } \,
     \dd \bcW_{t,j} ,
  \qquad t \in \RR_+ ,
 \end{equation}
 with initial value \ $\bcN_0 = \bzero$, \ and
 \begin{equation}\label{M_N}
  \bcM_{t,i}
  = \frac{1}{r}
    \int_0^t
     \sqrt{ \left[ \bm_\bxi^{r-i} \bPi 
                   \left(r \bcN_s
                         + s \sum_{\ell=1}^r \bm_\bxi^{\ell-1} \bm_\bvare\right)
                   \right]^+
            \odot \bV_\bxi } \,
     \dd \bcW_{s,i} , \qquad i \in \{1, \ldots, r\} .
 \end{equation}
If
 \[
   \bcN_t = \begin{bmatrix}
             \bcN_{t,1} \\
             \vdots \\
             \bcN_{t,r}
            \end{bmatrix} , \qquad
   \bcW_{t,j} = \begin{bmatrix}
                \bcW_{t,j,1} \\
                \vdots \\
                \bcW_{t,j,r}
               \end{bmatrix} , \qquad j \in \{1, \ldots, r\} , \qquad
   t \in \RR_+ ,
 \]
 denote the partitioning of \ $\bcN_t$ \ and \ $\bcW_{t,j}$,
 \ $j \in \{1, \ldots, r\}$, \ with respect to the partitioning
 \ $D_1$, \ldots, $D_r$ \ of the coordinates \ $\{1, \ldots, p\}$, \ then the
 process \ $(\bcN_t)_{t\in\RR_+}$ \ is the unique strong solution of the SDE
 \begin{equation}\label{SDE_N_mm}
  \dd \bcN_{t,i}
  = \sqrt{\left[\bv_i^\top
                \left(\bcN_{t,i}
                      + \frac{t}{r}
                        \sum_{\ell=i}^{i+r-1}
                         \tbm_{i,\ell} \bm_{\bvare,\ell}\right)\right]^+ }
    \sum_{\ell=i}^{i+r-1}
     \tbm_{i,\ell} \sqrt{(\tbm_{\ell+1,i} \bu_i) \odot \bV_{\ell+1}}
     \, \dd \bcW_{t,\ell+i-1,\ell+1} ,
 \end{equation}
 $i \in \{1, \ldots, r\}$, \ where the second and third subscripts of
 \ $\bcW_{t,\ell+i-1,\ell+1}$ \ are considered modulo \ $r$.
\end{Thm}

From \eqref{Mk} we obtain the recursion
 \begin{equation}\label{regr}
  \bX_{rk-i+1}
  = \bm_{\bxi}^r \bX_{rk-r-i+1}
    + \sum_{\ell=i}^r \bm_{\bxi}^{\ell-i} (\bM_{k,\ell} + \bm_{\bvare})
    + \sum_{\ell=1}^{i-1} \bm_{\bxi}^{\ell-i+r} (\bM_{k-1,\ell} + \bm_{\bvare})
 \end{equation} 
 for \ $k \in \NN$, \ $i \in \{1, \ldots, r\}$, \ where \ $\bM_{0,i} := \bzero$
 \ for \ $i \in \{1, \ldots, r\}$.
\ This recursion implies
 \begin{equation}\label{X}
  \bX_{rk-i+1}
  = \sum_{j=1}^k
     \bm_{\bxi}^{(k-j)r} 
     \left[ \sum_{\ell=i}^r \bm_{\bxi}^{\ell-i} (\bM_{j,\ell} + \bm_{\bvare})
            + \sum_{\ell=1}^{i-1}
               \bm_{\bxi}^{\ell-i+r} (\bM_{j-1,\ell} + \bm_{\bvare}) \right]
 \end{equation} 
 for \ $k \in \NN$, \ $i \in \{1, \ldots, r\}$.
\ Applying Lemma \ref{Conv2Funct}, which is a version of the continuous
 mapping theorem, together with \eqref{X}, \eqref{N} and Theorem \ref{Conv_M},
 we show the following convergence result.

\begin{Thm}\label{Conv_bX}
Under the assumptions of Theorem \ref{main}, we have
 \[
   \bcX^{(n)} \distr \bcX \qquad \text{as \ $n \to \infty$,}
 \]
 where
 \[
   \bcX_t
   = \begin{bmatrix}
      \bcX_{t,1} \\
      \bcX_{t,2} \\
      \vdots \\
      \bcX_{t,r}
     \end{bmatrix}
   = \begin{bmatrix}
      \bPi
      \sum_{\ell=1}^r
       \bm_\bxi^{\ell+r-1} (\bcM_{t,\ell} + r^{-1} t \bm_\bvare) \\
      \bPi
      \sum_{\ell=1}^r
       \bm_\bxi^{\ell+r-2} (\bcM_{t,\ell} + r^{-1} t \bm_\bvare) \\
      \vdots \\
      \bPi \sum_{\ell=1}^r \bm_\bxi^\ell (\bcM_{t,\ell} + r^{-1} t \bm_\bvare)
     \end{bmatrix}
   = \begin{bmatrix}
      \bm_\bxi^r
      \bPi \left( \bcN_t + \frac{t}{r}
                           \sum_{\ell=1}^r \bm_\bxi^{\ell-1} \bm_\bvare \right) \\
      \bm_\bxi^{r-1} \bPi
      \left( \bcN_t + \frac{t}{r}
                      \sum_{\ell=1}^r \bm_\bxi^{\ell-1} \bm_\bvare \right) \\
      \vdots \\
      \bm_\bxi
      \bPi \left( \bcN_t + \frac{t}{r}
                           \sum_{\ell=1}^r \bm_\bxi^{\ell-1} \bm_\bvare \right)
     \end{bmatrix}
 \]
 for all \ $t \in \RR_+$. 
\ Hence we obtain \ $\bcX_{t,i} = \bm_\bxi^{r-i+1} \bcY_t$,
 \ $i \in \{1, \ldots, r\}$, \ with
 \[
   \bcY_t := \bPi \left(\bcN_t
                        + \frac{t}{r}
                          \sum_{\ell=1}^r \bm_\bxi^{\ell-1} \bm_\bvare \right) ,
   \qquad t \in \RR_+ ,
 \]
 for which we have \ $\bm_\bxi^r \bcY_t = \bcY_t$, \ $t \in \RR_+$.
\end{Thm}

Theorem \ref{main} is an easy consequence of Theorems \ref{PROC_N} and
 \ref{Conv_bX}.
Indeed, \ $\bPi = \bPi_1 \oplus \cdots \oplus \bPi_r$ \ and
 \ $\bPi_i = r \bu_i \bv_i^\top$, \ $r \bv_i^\top \bu_i = 1$ \ for all
 \ $i \in \{1, \ldots, r\}$, \ hence we conclude from Theorems \ref{PROC_N}
 and \ref{Conv_bX} that for each \ $i \in \{1, \ldots, r\}$, \ the process
 \ $\cZ_{t,i} := \bv_i^\top \bcY_{t,i}$, \ $t \in \RR_+$, \ satisfies
 \[
   \cZ_{t,i}
   = \bv_i^\top \bPi_i
     \left(\bcN_{t,i} + \frac{t}{r}
                       \sum_{\ell=i}^{i+r-1} \tbm_{i,\ell} \bm_{\bvare,\ell} \right)
   = \bv_i^\top
     \left(\bcN_{t,i}
           + \frac{t}{r}
             \sum_{\ell=i}^{i+r-1} \tbm_{i,\ell} \bm_{\bvare,\ell} \right) ,
   \qquad t \in \RR_+ ,
 \]
 hence
 \[
   \cZ_{t,i} \bu_i
   = \bu_i \bv_i^\top
     \left(\bcN_{t,i}
           + \frac{t}{r} \sum_{\ell=i}^{i+r-1} \tbm_{i,\ell} \bm_{\bvare,\ell} \right)
   = \bPi_i
     \left(\bcN_{t,i} + \frac{t}{r}
                       \sum_{\ell=i}^{i+r-1} \tbm_{i,\ell} \bm_{\bvare,\ell} \right)
   = \bcY_{t,i} .
 \]
By It\^o's formula, we obtain that \ $(\cZ_{t,i})_{t\in\RR_+}$ \ is a strong
 solution of the SDE
 \begin{equation}\label{SDE_Z_bW}
  \begin{aligned}
   \dd \cZ_{t,i}
   &= r^{-1} \bv_i^\top \sum_{\ell=i}^{i+r-1} \tbm_{i,\ell} \bm_{\bvare,\ell} \, \dd t \\
   &\quad
      + \bv_i^\top \sqrt{r^{-1} \cZ_{t,i}^+}
        \sum_{\ell=i}^{i+r-1} 
         \tbm_{i,\ell}
         \sqrt{(\tbm_{\ell+1,i} \bu_i) \odot \bV_{\ell+1}}
         \, \dd \bcW_{t,\ell+i-1,\ell+1}
  \end{aligned}
 \end{equation}
 with initial value \ $\cZ_{0,i} = 0$. 
\ This equation can be written in the form \eqref{SDE_Z}, where
 \ $(\cW_{t,i})_{t\in\RR_+}$, \ $i \in \{1, \ldots, r\}$, \ are independent
 standard Wiener processes.
Indeed, we have
 \begin{multline*}
  \bv_i^\top
   \sum_{\ell=i}^{i+r-1}
    \tbm_{i,\ell} \left[ (\tbm_{\ell+1,i} \bu_i) \odot \bV_{\ell+1} \right]
    \tbm_{i,\ell}^\top \bv_i \\
  \begin{aligned}
   &= \sum_{\ell=i}^{i+r-1}
       \left(\bv_i^\top \tbm_{i,\ell}
             \sqrt{(\tbm_{\ell+1,i} \bu_i) \odot \bV_{\ell+1}} \right)
       \left(\bv_i^\top \tbm_{i,\ell}
             \sqrt{(\tbm_{\ell+1,i} \bu_i) \odot \bV_{\ell+1}} \right)^\top \\
   &= \sum_{\ell=i}^{i+r-1}
       \left\|\bv_i^\top \tbm_{i,\ell}
              \sqrt{(\tbm_{\ell+1,i} \bu_i) \odot \bV_{\ell+1}}\right\|^2
   \ne 0 .
  \end{aligned}
 \end{multline*}
Hence, if
 \ $\bv_i^\top \tbm_{i,\ell} \sqrt{(\tbm_{\ell+1,i} \bu_i) \odot \bV_{\ell+1}}
    = \bzero$ \ for each \ $\ell \in \{i, \ldots, i+r-1\}$, \ then
 \eqref{SDE_Z} trivially follows, and if there exists
 \ $\ell \in \{i, \ldots, i+r-1\}$ \ with
 \ $\bv_i^\top \tbm_{i,\ell} \sqrt{(\tbm_{\ell+1,i} \bu_i) \odot \bV_{\ell+1}}
    \ne \bzero$,
 \ then \eqref{SDE_Z} holds with
 \[
   \cW_{t,i}
   := \frac{\bv_i^\top
            \sum_{\ell=i}^{i+r-1}
             \tbm_{i,\ell} \sqrt{(\tbm_{\ell+1,i} \bu_i) \odot \bV_{\ell+1}}
             \, \bcW_{t,\ell+i-1,\ell+1}}
           {\bv_i^\top
            \sum_{\ell=i}^{i+r-1}
             \tbm_{i,\ell}
             \left[ (\tbm_{\ell+1,i} \bu_i) \odot \bV_{\ell+1} \right]
             \tbm_{i,\ell} ^\top \bv_i} ,
   \qquad t \in \RR_+ ,
   \qquad i \in \{1, \ldots, r\} ,
 \]
 which are independent standard Wiener processes, since
 \ $\bigl\{ (\ell+i-1,\ell+1) : \ell \in \{i, \ldots, i+r-1\} \bigr\}$,
 \ $i \in \{1, \ldots, r\}$, \ are disjoint sets.
Consequently, we conclude \eqref{Conv_X}.

\section{Proof of Theorem \ref{PROC_N}}
\label{Proof_Conv_M}

If \ $(\bcM_t)_{t\in\RR_+}$ \ is a strong solution of the SDE \eqref{SDE_M}, then
 the process \ $(\bcN_t)_{t\in\RR_+}$ \ is a strong solution of the SDE
 \eqref{SDE_N} with initial value \ $\bcN_0 = \bzero$, \ and \eqref{M_N}
 trivially holds.

Using the block matrix form of \ $\bm_\bxi$, \ $\bPi$ \ and
 \ $\bV_{\bxi_1}$, \ldots, $\bV_{\bxi_p}$ \ (see \eqref{partition_m}, \eqref{bPi}
 and \eqref{Cov_V}), we obtain
 \begin{equation}\label{SDE_N_m}
  \dd \bcN_{t,i}
  = \frac{1}{r}
    \sum_{j=1}^r
     \tbm_{i,i+j-1}
     \sqrt{ \left[ \tbm_{i-r+j,i} \bPi_i
                   \left(r \bcN_{t,i}
                         + t \sum_{\ell=i}^{i+r-1}
                              \tbm_{i,\ell} \bm_{\bvare,\ell} \right) \right]^+
            \odot \bV_{i-r+j} } \,
     \dd \bcW_{t,j,i+j}
 \end{equation}
 for each \ $i \in \{1, \ldots, r\}$.
\ Indeed, the covariance matrices \ $\bV_{\bxi_j}$ \ $j \in \{1, \ldots, p\}$,
 \ have block-diagonal form, see \eqref{Cov_V}, hence 
 \begin{multline*}
  \sqrt{ \left[ \bm_\bxi^{r-j} \bPi
                \left(r \bcN_t
                      + t \sum_{\ell=1}^r
                           \bm_\bxi^{\ell-1} \bm_\bvare\right) \right]^+
         \odot \bV_\bxi } \\
  = \bigoplus_{i=r-j+2}^{2r-j+1}
     \sqrt{ \left[ \tbm_{i-r+j,i} \bPi_i
                   \left(r \bcN_{t,i}
                         + t \sum_{\ell=i}^{r+i-1}
                              \tbm_{i,\ell} \bm_{\bvare,\ell} \right) \right]^+
            \odot \bV_{i-r+j} } ,
 \end{multline*}
 where we also used that for an arbitrary matrix \ $\bA \in \RR^{p \times p}$
 \ with partitioning
 \[
   \bA = \begin{bmatrix}
          \bA_{1,1} & \cdots & \bA_{1,r} \\
          \vdots & \ddots & \vdots \\
          \bA_{r,1} & \cdots & \bA_{r,r}
         \end{bmatrix}
 \]
 with respect to the partitioning \ $D_1$, \ldots, $D_r$ \ of the coordinates
 \ $\{1, \ldots, p\}$, \ we have
 \[
   \bm_\bxi^k \bA
   = \begin{bmatrix}
      \tbm_{1,k+1} \bA_{k+1,1} & \cdots & \tbm_{1,k+1} \bA_{k+1,r} \\
      \vdots & \ddots & \vdots \\
      \tbm_{r,k+r} \bA_{k+r,1} & \cdots & \tbm_{r,k+r} \bA_{k+r,r}
     \end{bmatrix} 
 \]
 for all \ $k \in \{1, \ldots r-1\}$, \ where the subscripts are considered
 modulo \ $r$.
\ Substituting this into \eqref{SDE_N} and using again the above block form of
 \ $\bm_\bxi^k \bA$ \ for \ $\bA \in \RR^{p \times p}$ \ and
 \ $k \in \{1, \ldots r-1\}$, \ we obtain \eqref{SDE_N_m}.
Using \ $\bPi_i = r \bu_i \bv_i^\top$ \ for all \ $i \in \{1, \ldots, r\}$,
 \ equation \eqref{SDE_N_m} can be written in the form \eqref{SDE_N_mm}.
\proofend

\section{Proof of Theorem \ref{Conv_M}}
\label{Proof_Conv_M}

In order to prove \ $\bcM^{(n)} \distr \bcM$, \ we want to apply Theorem
 \ref{Conv2DiffThm} for \ $\bcU = \bcM$, \ $\bU^{(n)}_k = n^{-1} \bM_k$ \ and
 \ $\cF^{(n)}_k := \cF_k$ \ for \ $n \in \NN$ \ and \ $k \in \ZZ_+$, \ and with
 coefficient function
 \ $\gamma : \RR_+ \times (\RR^p)^r \to (\RR^{p \times p})^{r \times r}$
 \ of the SDE \eqref{SDE_M} given by
 \[
   \gamma(t, \bx)
   = \frac{1}{r}
     \bigoplus_{i=1}^r
      \sqrt{ \left[ \bm_\bxi^{r-i} \bPi
                    \sum_{j=1}^r
                     \bm_\bxi^{j-1} (r \bx_j + t \bm_\bvare) \right]^+
             \odot \bV_\bxi } , \qquad
   \bx = \begin{bmatrix} \bx_1 \\ \vdots \\ \bx_r \end{bmatrix} \in (\RR^p)^r .
 \] 
The aim of the following discussion is to show that the SDE \eqref{SDE_M} has
 a unique strong solution \ $\bigl(\bcM_t^{(\bx_0)}\bigr)_{t\in\RR_+}$ \ with
 initial value \ $\bcM_0^{(\bx_0)} = \bx_0$ \ for all \ $\bx_0 \in (\RR^p)^r$. 
\ Clearly, it is sufficient to prove that the SDE \eqref{SDE_N} has
 a unique strong solution \ $(\bcN_t^{(\by_0)})_{t\in\RR_+}$ \ with initial value
 \ $\bcN_0^{(\by_0)} = \by_0$ \ for all \ $\by_0 \in (\RR^p)^r$.
\ Indeed, if \ $\bigl(\bcM_t^{(\bx_0)}\bigr)_{t\in\RR_+}$ \ is a strong solution of
 the SDE \eqref{SDE_M} with initial value \ $\bcM_0^{(\bx_0)} = \bx_0$ \ with
 some \ $\bx_0 \in (\RR^p)^r$ \ then
 \ $\bcN_t := \sum_{i=1}^r \bm_\bxi^{i-1} \bcM_{t,i}^{(\bx_0)}$ \ is a strong
 solution of the SDE \eqref{SDE_N} with initial value
 \ $\sum_{i=1}^r \bm_\bxi^{i-1} \bx_{0,i}$. 
\ Conversely, if \ $\bigl(\bcN_t^{(\by_0)}\bigr)_{t\in\RR_+}$ \ is a strong
 solution of the SDE \eqref{SDE_N} with initial value
 \ $\bcN_0^{(\by_0)} = \by_0$ \ with some \ $\by_0 \in \RR^p$ \ then there exists
 \ $\bx_0 \in (\RR^p)^r$ \ such that
 \ $\by_0 = \sum_{i=1}^r \bm_\bxi^{i-1} \bx_{0,i} \in \RR^p$ \ (for instance,
 \ $\bx_{0,1} = \by_0$ \ and \ $\bx_{0,2} = \ldots = \bx_{0,r} = \bzero$), \ and
 \[
   \bcM_{t,i}
    := \bx_{0,i}
       + \frac{1}{r}
         \int_0^t
          \sqrt{\left[ \bm_\bxi^{r-i} \bPi
                       \left(r \bcN_s^{(\by_0)}
                             + s \sum_{\ell=1}^r \bm_\bxi^{\ell-1} \bm_\bvare\right)
                \right]^+
                \odot \bV_\bxi}
          \, \dd \bcW_{s,i} , \qquad i \in \{1, \ldots, r\} ,
 \]
 is a strong solution of the SDE \eqref{SDE_M} with initial value \ $\bx_0$.

Hence it is enough to show that the SDE \eqref{SDE_N} has a unique strong
 solution \ $\bigl(\bcN_t^{(\by_0)}\bigr)_{t\in\RR_+}$ \ with
 initial value \ $\bcN_0^{(\by_0)} = \by_0$ \ for all \ $\by_0 \in \RR^p$.
\ First observe that if \ $\bigl(\bcN_{t,i}^{(\by_{0,i})}\bigr)_{t\in\RR_+}$ \ is a
 strong solution of the SDE \eqref{SDE_N_mm} with initial value
 \ $\bcN_{0,i}^{(\by_{0,i})} = \by_{0,i} \in \RR^p$, \ then, by It\^o's formula, the
 process \ $(\cP_{t,i}, \, \bcQ_{t,i})_{t\in\RR_+}$, \ defined by
 \[
   \cP_{t,i} := \bv_i^\top
               \left( \bcN_{t,i}^{(\by_{0,i})}
                      + \frac{t}{r}
                        \sum_{\ell=i}^{i+r-1}
                         \tbm_{i,\ell} \bm_{\bvare,\ell} \right) , \qquad
   \bcQ_{t,i} := \bcN_{t,i}^{(\by_{0,i})} - \cP_{t,i} \bu_i
 \]
 is a strong solution of the SDE
 \begin{equation}\label{SDE_P_Q}
  \begin{cases}
   \dd \cP_{t,i}
   &= \frac{1}{r}
      \bv_i^\top \sum_{\ell=i}^{i+r-1} \tbm_{i,\ell} \bm_{\bvare,\ell} \, \dd t \\
   &\quad
      + \sqrt{r^{-1} \cP_{t,i}^+} \, \bv_i^\top
        \sum_{\ell=i}^{i+r-1}
         \tbm_{i,\ell} \sqrt{(\tbm_{\ell+1,i+r} \bu_i) \odot \bV_{\ell+1}} 
         \, \dd \bcW_{t,\ell+i-1,\ell+1} , \\[2mm]
   \dd \bcQ_{t,i}
   &= - \frac{1}{r} \bPi_i
        \sum_{\ell=i}^{i+r-1} \tbm_{i,\ell} \bm_{\bvare,\ell} \, \dd t \\
   &\quad
      + \sqrt{r^{-1} \cP_{t,i}^+} \, (\bI_p - \bPi_i)
        \sum_{\ell=i}^{i+r-1}
         \tbm_{i,\ell} \sqrt{(\tbm_{\ell+1,i+r} \bu_i) \odot \bV_{\ell+1}} 
         \, \dd \bcW_{t,\ell+i-1,\ell+1}
  \end{cases}
 \end{equation} 
 with initial value
 \ $(\cP_{0,i}, \, \bcQ_{0,i})
    = (\bv_i^\top \by_{0,i}, (\bI_p - \bPi_i) \by_{0,i})$,
 \ where \ $\bI_p \in \RR^{p \times p}$ \ denotes the unit matrix.
Indeed, the first SDE of \eqref{SDE_P_Q} is an easy consequence of the SDE
 \eqref{SDE_N_mm}.
The second one can be checked as follows. 
By It\^o's formula,
 \begin{align*}
  \dd \bcQ_{t,i}
  &= \dd \bcN_{t,i}^{(\by_{0,i})} - \bu_i \, \dd \cP_{t,i}
   = \dd \bcN_{t,i}^{(\by_{0,i})}
     - \bu_i \bv_i^\top
       \left( \dd \bcN_{t,i}^{(\by_{0,i})}
              + \frac{1}{r}
                \sum_{\ell=i}^{i+r-1}
                 \tbm_{i,\ell} \bm_{\bvare,\ell} \, \dd t \right) \\
  &= - \frac{1}{r} \bPi_i
       \sum_{\ell=i}^{i+r-1} \tbm_{i,\ell} \bm_{\bvare,\ell} \, \dd t
     + (\bI_p - \bPi_i) \, \dd \bcN_{t,i}^{(\by_{0,i})}
 \end{align*}
 with
 \ $\bcQ_{0,i} = \bcN_{0,i}^{(\by_{0,i})} - \cP_{0,i} \bu_i
    = \by_{0,i} - (\bv_i^\top \by_{0,i}) \bu_i
    = \by_{0,i} - \bu_i \bv_i^\top \by_{0,i} = (\bI_p - \bPi_i) \by_{0,i}$.
\ Conversely, if 
 \ $(\cP_{t,i}^{(p_{0,i}, \bq_{0,i})}, \, \bcQ_{t,i}^{(p_{0,i}, \bq_{0,i})})_{t\in\RR_+}$
 \ is a strong solution of the SDE \eqref{SDE_P_Q} with initial value
 \ $\bigl(\cP_{0,i}^{(p_{0,i}, \bq_{0,i})}, \, \bcQ_{0,i}^{(p_{0,i}, \bq_{0,i})}\bigr)
    = (p_{0,i}, \bq_{0,i}) \in \RR \times \RR^p$,
 \ then, again by It\^o's formula,
 \[
   \bcN_{t,i}
   := \cP_{t,i}^{(p_{0,i}, \bq_{0,i})} \, \bu_i + \bcQ_{t,i}^{(p_{0,i}, \bq_{0,i})} ,
   \qquad t \in \RR_+ , 
 \] 
 is a strong solution of the SDE \eqref{SDE_N_mm} with initial value
 \ $\bcN_{0,i} = p_{0,i} \bu_i + \bq_{0,i}$.
\ The correspondence
 \ $\by_{0,i} \leftrightarrow (p_{0,i}, \bq_{0,i})
    := (\bv_i^\top \by_{0,i}, \, (\bI_p - \bPi_i) \by_{0,i})$
 \ is a bijection between \ $\RR^p$ \ and
 \ $\RR \times \{ \bq \in \RR^p : \bv_i^\top \bq = 0 \}$, \ since
 \ $\by_{0,i} = p_{0,i} \bu_i + \bq_{0,i}$, \ and
 \[
   (\bI_p - \bPi_i) (p_{0,i} \bu_i + \bq_{0,i})
   = p_{0,i} \bu_i + \bq_{0,i} - \bPi_i p_{0,i} \bu_i + \bPi_i \bq_{0,i}
   = p_{0,i} \bu_i + \bq_{0,i} - p_{0,i} \bu_i \bv_i^\top \bu_i
     + \bu_i \bv_i^\top \bq_{0,i}
   = \bq_{0,i} .
 \]
The right hand side of the SDE \eqref{SDE_P_Q} contains only the process
 \ $(\cP_{t,i})_{t\in\RR_+}$, \ hence it is enough to show that the first equation
 of \eqref{SDE_P_Q} has a unique strong solution
 \ $(\cP_{t,i}^{(p_{0,i}, \bq_{0,i})})_{t\in\RR_+}$ \ with initial value
 \ $\cP_{0,i}^{(p_{0,i}, \bq_{0,i})} = p_{0,i}$ \ for all \ $p_{0,i} \in \RR$. 
\ The first equation of \eqref{SDE_P_Q} is the same as \eqref{SDE_Z_bW}, which
 can be written in the form \eqref{SDE_Z}, see the end of Section
 \ref{Convergence_result}.
Hence, by Remark \ref{Rem_SDE_Z}, the first equation of the SDE \eqref{SDE_P_Q}
 has a unique strong solution \ $(\cP_{t,i}^{(p_{0,i})})_{t\in\RR_+}$ \ with initial
 value \ $\cP_{0,i}^{(p_{0,i})} = p_{0,i}$ \ for all \ $p_{0,i} \in \RR$.
\ Consequently, the SDE \eqref{SDE_P_Q}, and hence the SDE \eqref{SDE_M} admit
 a unique strong solution with arbitrary initial value.

Now we show that conditions (i) and (ii) of Theorem \ref{Conv2DiffThm} hold.
We have to check that for each \ $T > 0$,
 \begin{align} \label{Cond1}
  &\sup_{t\in[0,T]}
    \bigg\| \frac{1}{(nr)^2} \sum_{k=1}^{\nt}
            \EE\bigl[ \bM_k \bM_k^\top \, \big| \, \cF_{rk-r} \bigr]
            - \int_0^t (\bcR^{(n)}_s)^+ \, \dd s \bigg\|
   \stoch 0,\\
  &\frac{1}{(nr)^2}
   \sum_{k=1}^{\nT}
    \EE\bigl( \|\bM_k\|^2 \bbone_{\{\|\bM_k\| > n \theta\}} \,
              \big| \, \cF_{rk-r} \bigr)
   \stoch 0   \qquad\text{for all \ $\theta>0$} \label{Cond2}
 \end{align}
 as \ $n \to \infty$, \ where the process \ $(\bcR^{(n)}_s)_{s \in \RR_+}$ \ is
 defined by
 \begin{equation}\label{Rnt}
  \bcR^{(n)}_s
  := \frac{1}{r^2}
     \bigoplus_{i=1}^r
      \left\{ \left[ \bm_\bxi^{r-i} \bPi
                     \sum_{j=1}^r
                      \bm_\bxi^{j-1}
                      (\bcM_{s,j}^{(n)} + r^{-1} s \bm_\bvare) \right]
              \odot \bV_\bxi \right\}
 \end{equation}
 for \ $s \in \RR_+$, \ $n \in \NN$.
\ By \eqref{Mk},
 \begin{align*}
  &\bPi \sum_{j=1}^r \bm_\bxi^{j-1} (\bcM_{s,j}^{(n)} + r^{-1} s \bm_\bvare) \\
  &\qquad\qquad
   = \bPi
       \sum_{j=1}^r
        \bm_\bxi^{j-1} 
        \Biggl((nr)^{-1} \sum_{k=1}^\ns
                         (\bX_{rk-j+1} - \bm_\bxi \bX_{rk-j} - \bm_\bvare)
                + r^{-1} s \bm_\bvare \Biggr) \\
  &\qquad\qquad
   = (nr)^{-1} \bPi
       \sum_{k=1}^\ns
        \sum_{j=1}^r
         (\bm_\bxi^{j-1} \bX_{rk-j+1} - \bm_\bxi^j \bX_{rk-j}
          - \bm_\bxi^{j-1} \bm_\bvare)
       + r^{-1} s \bPi \sum_{j=1}^r \bm_\bxi^{j-1} \bm_\bvare \\
  &\qquad\qquad
   = (nr)^{-1} \bPi
       \sum_{k=1}^\ns
        \left(\bX_{rk} - \bm_\bxi^r \bX_{rk-r}
              - \sum_{j=1}^r \bm_\bxi^{j-1} \bm_\bvare\right)
       + r^{-1} s \bPi \sum_{j=1}^r \bm_\bxi^{j-1} \bm_\bvare
 \end{align*}
 \begin{align*}
  &\qquad\qquad
   = (nr)^{-1}
       \sum_{k=1}^\ns
        \Biggl(\bPi \bX_{rk} - \bPi \bX_{rk-r}
               - \bPi \sum_{j=1}^r \bm_\bxi^{j-1} \bm_\bvare\Biggr)
       + r^{-1} s \bPi \sum_{j=1}^r \bm_\bxi^{j-1} \bm_\bvare \\
  &\qquad\qquad
   = (nr)^{-1} \bPi \bX_{r\ns}
       + \left( s - \frac{\ns}{n} \right) r^{-1} \bPi
         \sum_{j=1}^r \bm_\bxi^{j-1} \bm_\bvare ,
 \end{align*}
 where we used
 \begin{equation}\label{Pimr}
   \bPi \bm_\bxi^r
   = \left( \lim\limits_{n \to \infty} \bm_\bxi^{nr} \right) \bm_\bxi^r
   = \lim\limits_{n \to \infty} \bm_\bxi^{(n+1)r} = \bPi .
 \end{equation}
Consequently,
 \[
   \bcR^{(n)}_s
   = \frac{1}{r^2}
     \bigoplus_{i=1}^r
      \left\{ \left[ n^{-1} \bm_\bxi^{r-i} \bPi \bX_{r\ns}
                     + \left( s - \frac{\ns}{n} \right)
                       \bPi \sum_{j=1}^r \bm_\bxi^{j-1} \bm_\bvare \right]
              \odot \bV_\bxi \right\} ,
 \]
 since
 \[
   \bm_\bxi^{r-i} \bPi
    = \bm_\bxi^{r-i} \left( \lim\limits_{n \to \infty} \bm_\bxi^{nr} \right)
    = \left( \lim\limits_{n \to \infty} \bm_\bxi^{nr+r-i} \right)
    = \left( \lim\limits_{n \to \infty} \bm_\bxi^{nr} \right) \bm_\bxi^{r-i}
    = \bPi \bm_\bxi^{r-i}
 \]
 and \eqref{Pimr} implies
 \begin{align*}
  \bm_\bxi^{r-i} \bPi \sum_{j=1}^r \bm_\bxi^{j-1}
  &= \bPi \bm_\bxi^{r-i} \sum_{j=1}^r \bm_\bxi^{j-1} \\
  &= \bPi \left( \sum_{j=1}^i \bm_\bxi^{j-1+r-i}
                 + \bm_\bxi^r \sum_{j=i+1}^r \bm_\bxi^{j-1-i} \right)
   = \bPi \sum_{j=1}^r \bm_\bxi^{j-1} .
 \end{align*}
Thus \ $(\bcR^{(n)}_t)^+ = \bcR^{(n)}_t$, \ and 
 \begin{align*}
  \int_0^t (\bcR^{(n)}_s)^+ \,\dd s
  = \bigoplus_{i=1}^r
     \bigg\{ \bigg[ & \frac{1}{(nr)^2} \bm_\bxi^{r-i} \bPi
                      \sum_{\ell=0}^{\nt-1} \bX_{r \ell}
                      + \frac{nt-\nt}{(nr)^2}
                        \bm_\bxi^{r-i} \bPi \bX_{r \nt} \\
                    & + \frac{\nt+(nt-\nt)^2}{2 (nr)^2} \bPi
                        \sum_{j=1}^r \bm_\bxi^{j-1} \bm_\bvare \bigg]
             \odot \bV_\bxi \bigg\} .
 \end{align*}
Using \eqref{Mcond}, we obtain
 \[
   \frac{1}{(nr)^2} \sum_{k=1}^{\nt} \EE( \bM_k \bM_k^\top \mid \cF_{rk-r} )
   = \bigoplus_{i=1}^r
      \left\{ \frac{\nt}{n^2} \bV_\bvare
              + \frac{1}{n^2}
                \sum_{k=1}^{\nt}
                 \bigl[ \bm_\bxi^{r-i} \bX_{rk-r}
                        + \sum_{j=1}^{r-i} \bm_\bxi^{j-1} \bm_\bvare \bigr]
                 \odot \bV_\bxi \right\} .
 \]
Hence, taking into acount that \ $\bX_0 = \bzero$, \ in order to show
 \eqref{Cond1}, it suffices to prove
 \begin{equation}\label{Cond11}
  \frac{1}{n^2}
  \sup_{t \in [0,T]}
   \sum_{k=1}^\nt
    \| (\bI_p - \bPi) \bX_{rk} \|
  \stoch 0,\qquad
  \frac{1}{n^2}
  \sup_{t \in [0,T]}
   \|\bX_{r \nt}\|
  \stoch 0
 \end{equation}
 as \ $n \to \infty$.
\ Using \eqref{X} and \eqref{Pimr}, \ we obtain
 \begin{align*} 
  (\bI_p - \bPi) \bX_{rk}
  &= (\bI_p - \bPi)
     \sum_{j=1}^k
      \bm_\bxi^{(k-j)r}
      \sum_{i=1}^r \bm_\bxi^{i-1} (\bM_{j,i} + \bm_\bvare) \\
  &= \sum_{j=1}^k
      \left( \bm_\bxi^{(k-j)r} - \bPi \right)
      \sum_{i=1}^r \bm_\bxi^{i-1} (\bM_{j,i} + \bm_\bvare) .
 \end{align*}
Hence by \eqref{rate},
 \begin{align*} 
  \sum_{k=1}^{\nt}
   \| (\bI_p - \bPi) \bX_{rk} \|
   & \leq c
          \sum_{k=1}^\nt
           \sum_{j=1}^k
            \kappa^{k-j}
            \sum_{i=1}^r
             \|\bm_\bxi^{i-1}\| \left\| \bM_{j,i} + \bm_\bvare \right\| \\
   & \leq c K
          \sum_{j=1}^\nt
           \sum_{k=j}^\nt
            \kappa^{k-j}
            \left( \sum_{i=1}^r \|\bM_{j,i}\| + r \| \bm_\bvare \| \right) \\
   & \leq \frac{c K}{1 - \kappa}
          \left( \sum_{j=1}^{\nt} \sum_{i=1}^r \|\bM_{j,i}\|
                 + r \nt \| \bm_\bvare \| \right) ,
 \end{align*}
 where \ $K := \max_{i\in\{1,\ldots,r\}} \|\bm_\bxi^{i-1}\|$.
\ Moreover, by \eqref{X},
 \begin{align*}
  \|\bX_{r \nt}\|
  & \leq \sum_{j=1}^{\nt}
          \| \bm_\bxi^{(\nt-j)r} \| \,
          \sum_{i=1}^r \| \bm_\bxi^{i-1} \|  \| \bM_{j,i} + \bm_\bvare \| \\
  & \leq K (c + \|\bPi\|)
         \left( r \nt \| \bm_\bvare \|
                + \sum_{j=1}^{\nt} \sum_{i=1}^r \|\bM_{j,i}\| \right) ,
 \end{align*}
 since
 \ $\| \bm_\bxi^{(\nt-j)r} \| \leq \| \bm_\bxi^{(\nt-j)r} - \bPi \| + \|\bPi\|
    \leq c + \|\bPi\|$
 \ by \eqref{rate}.
Consequently, in order to prove \eqref{Cond11}, it suffices to show
 \[
   \frac{1}{n^2} \sum_{j=1}^{\nT} \|\bM_{j,i}\| \stoch 0 \qquad
   \text{as \ $n \to \infty$ \ for all \ $i \in \{1, \ldots, r\}$.}
 \]
In fact, Lemma \ref{EEX} yields
 \ $n^{-2} \sum_{j=1}^{\nT} \EE(\|\bM_{j,i}\|) \to 0$, \ $i \in \{1, \ldots, r\}$,
 \ thus we obtain \eqref{Cond1}.

Next we check condition \eqref{Cond2}.
We have
 \[
   \EE\bigl( \|\bM_k\|^2 \bbone_{\{\|\bM_k\| > n \theta\}} \,
              \big| \, \cF_{k-1} \bigr)
   \leq n^{-2} \theta^{-2}
        \EE\bigl( \|\bM_k\|^4 \, \big| \, \cF_{k-1} \bigr) ,
 \]  
 hence \ $n^{-4} \sum_{k=1}^{\nT} \EE\bigl(\|\bM_k\|^4\bigr) \to 0$ \ as
 \ $n \to \infty$, \ since
 \ $\EE\bigl(\|\bM_k\|^4\bigr) = \OO(k^2)$ \ by Lemma \ref{EEX}.
This yields \eqref{Cond2}.
\proofend

\section{Proof of Theorem \ref{Conv_bX}}
\label{Proof_Conv_bX}

In order to prove Theorem \ref{Conv_bX}, we want to apply Lemma
 \ref{Conv2Funct} using Theorem \ref{Conv_M}.
By \eqref{X}, \ $\bcX^{(n)} = \Psi^{(n)}(\bcM^{(n)})$, \ where the mapping
 \[
   \Psi^{(n)} = \begin{bmatrix}
                \Psi^{(n)}_1 \\
                \vdots \\
                \Psi^{(n)}_r \\
               \end{bmatrix}
   : \DD(\RR_+, (\RR^p)^r) \to \DD(\RR_+, (\RR^p)^r)
 \]
 is given by
 \begin{align*}
  \Psi^{(n)}_i(f)(t)
  &:= \bm_{\bxi}^{\nt r}
      \Bigg[ \sum_{\ell=i}^r \bm_{\bxi}^{\ell-i} f_\ell(0)
             + \sum_{\ell=1}^{i-1} \bm_{\bxi}^{\ell-i+r} f_\ell(0) \Bigg] \\
  &\quad\:
      + \sum_{j=1}^\nt
         \bm_{\bxi}^{(\nt-j)r} 
         \Bigg[ \sum_{\ell=i}^r
                  \bm_{\bxi}^{\ell-i}
                  \left( f_\ell\left(\frac{j}{n}\right)
                         - f_\ell\left(\frac{j-1}{n}\right)
                         + \frac{1}{nr} \bm_{\bvare} \right) \\
  &\phantom{\quad\: + \sum_{j=1}^\nt \bm_{\bxi}^{(\nt-j)r} \Bigg[}
      + \sum_{\ell=1}^{i-1}
         \bm_{\bxi}^{\ell-i+r}
         \left( f_\ell\left(\frac{j-1}{n}\right)
                - f_\ell\left(\frac{j-2}{n}\right)
                + \frac{1}{nr} \bm_{\bvare} \right) \Bigg]
 \end{align*}
 for
 \[
   f = \begin{bmatrix}
        f_1 \\
        \vdots \\
        f_r \\
       \end{bmatrix} \in \DD(\RR_+, (\RR^p)^r) ,
   \qquad t \in \RR_+ , \qquad  n \in \NN , \qquad i \in \{1, \ldots, r\} .
 \]
Further, \ $\bcX = \Psi(\bcM)$, \ where the mapping
 \ $\Psi : \DD(\RR_+, (\RR^p)^r) \to \DD(\RR_+, (\RR^p)^r)$ \ is given by
 \[
   \Psi(f)(t)
   := \begin{bmatrix}
       \bPi \sum_{\ell=1}^r
             \bm_\bxi^{\ell+r-1} ( f_\ell(t) + r^{-1} t \bm_\bvare ) \\
       \bPi \sum_{\ell=1}^r
             \bm_\bxi^{\ell+r-2} ( f_\ell(t) + r^{-1} t \bm_\bvare ) \\
       \vdots \\
       \bPi \sum_{\ell=1}^r \bm_\bxi^\ell ( f_\ell(t) + r^{-1} t \bm_\bvare )
      \end{bmatrix} , \qquad
   f \in \DD(\RR_+, (\RR^p)^r), \qquad t \in \RR_+ .
 \]
Measurability of the mappings \ $\Psi^{(n)}$, \ $n \in \NN$, \ and \ $\Psi$
 \ can be checked as in Barczy et al. \cite[page 603]{BarIspPap0}, see Lemma
 \ref{measurability}.

The aim of the following discussion is to show that the set
 \[
   C := \Bigg\{ f \in \CC(\RR_+, (\RR^p)^r)
               : \text{$(\bI_{rp} - \bPi)
                        \Bigg[ \sum_{\ell=i}^r \bm_{\bxi}^{\ell-i} f_\ell(0)
                               + \sum_{\ell=1}^{i-1}
                                  \bm_{\bxi}^{\ell-i+r} f_\ell(0) \Bigg] = \bzero$
                        \ for all $i \in \{1, \ldots, r\}$} \Bigg\}
 \]
 satisfies \ $C \in \cD_\infty(\RR_+, (\RR^p)^r)$,
 \ $C \subset C_{\Psi, \, (\Psi^{(n)})_{n \in \NN}}$ \ and \ $\PP(\bcM \in C) = 1$.

First note that
 \ $C = \CC(\RR_+, (\RR^p)^r)
        \cap \pi_0^{-1}\Bigl(\bigl(\bI_{rp}
                                  - \bPi^{\oplus r}
                                  \mathrm{circ}_r(\bm_\bxi))^{-1}(\{ \bzero \}
                            \bigr)
                     \Bigr)$,
 \ where \ $\mathrm{circ}_r(\bm_\bxi)$ \ denotes the circulant matrix
 \[
   \mathrm{circ}_r(\bm_\bxi)
   := \begin{bmatrix}
       \bI_p & \bm_\bxi & \cdots & \bm_\bxi^{r-2} & \bm_\bxi^{r-1} \\
       \bm_\bxi^{r-1} & \bI_p & \cdots & \bm_\bxi^{r-3} & \bm_\bxi^{r-2} \\
       \vdots & \vdots & \ddots & \vdots & \vdots \\
       \bm_\bxi^2 & \bm_\bxi^3 & \cdots & \bI_p & \bm_\bxi \\
       \bm_\bxi & \bm_\bxi^2 & \cdots & \bm_\bxi^{r-1} & \bI_p
      \end{bmatrix} ,
 \]
 $\bPi^{\oplus r} := \bPi \oplus \cdots \oplus \bPi$, \ and
 \ $\pi_0 : \DD(\RR_+, (\RR^p)^r) \to (\RR^p)^r$ \ denotes the projection
 defined by \ $\pi_0(f) := f(0)$ \ for \ $f \in \DD(\RR_+, (\RR^p)^r)$.
\ Using that \ $\CC(\RR_+, (\RR^p)^r) \in \cD_\infty(\RR_+, (\RR^p)^r)$ \ (see,
 e.g., Ethier and Kurtz \cite[Problem 3.11.25]{EK} and Lemma \ref{closed}),
 the linear mapping
 \ $(\RR^p)^r \ni \bx
    \mapsto \bigl(\bI_{rp} - \bPi^{\oplus r} \mathrm{circ}_r(\bm_\bxi)\bigr) \bx
    \in (\RR^p)^r$
 \ is measurable and that \ $\pi_0$ \ is measurable (see, e.g., Ethier and
 Kurtz \cite[Proposition 3.7.1]{EK}), we obtain
 \ $C \in \cD_\infty(\RR_+, (\RR^p)^r)$.
 
Fix a function \ $f \in C$ \ and a sequence \ $(f^{(n)})_{n\in\NN}$ \ in
 \ $\DD(\RR_+, (\RR^p)^r)$ \ with \ $f^{(n)} \lu f$.
\ By the definition of \ $\Psi$, \ we have
 \ $\Psi(f) \in \CC(\RR_+, (\RR^p)^r)$.
\ Further, we can write
 \begin{multline*}
  \Psi^{(n)}_i(f^{(n)})(t)
  =\bPi
   \sum_{\ell=i}^r
    \bm_\bxi^{\ell-i}
    \left[ f^{(n)}_\ell\left(\frac{\nt}{n}\right)
           + \frac{\nt}{nr} \bm_\bvare \right] \\
  \begin{aligned}
   &+\bPi
     \sum_{\ell=1}^{i-1}
      \bm_\bxi^{\ell-i+r}
      \left[ f^{(n)}_\ell\left(\frac{\nt-1}{n}\right)
             + \frac{\nt}{nr} \bm_\bvare \right] \\
   &+\left(\bm_\bxi^{\nt r} -\bPi\right)
     \Bigg[ \sum_{\ell=i}^r \bm_{\bxi}^{\ell-i} f^{(n)}_\ell(0)
            + \sum_{\ell=1}^{i-1} \bm_{\bxi}^{\ell-i+r} f^{(n)}_\ell(0) \Bigg] \\
   &+\sum_{j=1}^{\nt}
      \bigl( \bm_\bxi^{(\nt-j)r} - \bPi \bigr)
      \Bigg[ \sum_{\ell=i}^r
              \bm_{\bxi}^{\ell-i}
              \left( f^{(n)}_\ell\left(\frac{j}{n}\right)
                     - f^{(n)}_\ell\left(\frac{j-1}{n}\right)
                     + \frac{1}{nr} \bm_{\bvare} \right) \\
   &\phantom{+ \sum_{j=2}^{\nt} \bigl( \bm_\bxi^{(\nt-j)r} - \bPi \bigr) \Bigg[}
            + \sum_{\ell=1}^{i-1}
               \bm_{\bxi}^{\ell-i+r}
               \left( f^{(n)}_\ell\left(\frac{j-1}{n}\right)
                      - f^{(n)}_\ell\left(\frac{j-2}{n}\right)
                      + \frac{1}{nr} \bm_{\bvare} \right) \Bigg] ,
  \end{aligned}
 \end{multline*}
 hence we have
 \begin{align*}
  &\|\Psi_i^{(n)}(f^{(n)})(t)-\Psi_i(f)(t)\|
   \leq K \| \bPi \|
        \sum_{\ell=1}^{r-i+1}
         \left( \left\| f_\ell^{(n)}\left(\frac{\nt}{n}\right)
                        - f_\ell(t) \right\|
                + \frac{1}{nr} \|\bm_{\bvare}\| \right) \\
  &\qquad\qquad
   + K \| \bPi \|
           \sum_{\ell=1}^{r-i+1}
            \left( \left\| f_\ell^{(n)}\left(\frac{\nt-1}{n}\right)
                           - f_\ell(t) \right\|
                   + \frac{1}{nr} \|\bm_{\bvare}\| \right) \\
  &\qquad\qquad
   + \left\| \left(\bm_\bxi^{\nt r} -\bPi\right)
              \Bigg[ \sum_{\ell=i}^r \bm_{\bxi}^{\ell-i} f^{(n)}_\ell(0)
                     + \sum_{\ell=1}^{i-1}
                        \bm_{\bxi}^{\ell-i+r} f^{(n)}_\ell(0) \Bigg] \right\|
 \end{align*}
 \begin{align*}
  &\qquad\qquad
   + K \sum_{j=1}^{\nt}
         \bigl\| \bm_\bxi^{(\nt-j)r} - \bPi \bigr\|
         \Bigg[ \sum_{\ell=i}^r
                 \left( \left\| f_\ell^{(n)}\left(\frac{j}{n}\right)
                                - f_\ell^{(n)}\left(\frac{j-1}{n}\right) \right\|
                        + \frac{1}{nr} \|\bm_{\bvare}\| \right) \\
  &\phantom{\qquad\qquad
            + K \sum_{j=1}^{\nt} \bigl\| \bm_\bxi^{(\nt-j)r} - \bPi \bigr\|
             \Bigg[\:}
                + \sum_{\ell=1}^{i-1}
                   \left( \left\| f_\ell^{(n)}\left(\frac{j-1}{n}\right)
                                  - f_\ell^{(n)}\left(\frac{j-2}{n}\right)
                          \right\|
                          + \frac{1}{nr} \|\bm_{\bvare}\| \right) \Bigg] .
 \end{align*}
Here for all \ $T > 0$ \ and \ $t \in [0, T]$,
 \begin{align*}
  \left\| f^{(n)}\left(\frac{\nt}{n}\right) - f(t) \right \|
  &\leq \left\| f^{(n)}\left(\frac{\nt}{n}\right)
                - f\left(\frac{\nt}{n}\right) \right\|
        + \left\| f\left(\frac{\nt}{n}\right) - f(t) \right\| \\
  &\leq \omega_T(f, n^{-1}) + \sup_{t \in [0,T]} \|f^{(n)}(t) - f(t)\| ,
 \end{align*}
 where \ $\omega_T(f, \cdot)$ \ is the modulus of continuity of \ $f$ \ on
 \ $[0, T]$, \ and we have \ $\omega_T(f, n^{-1}) \to 0$ \ since \ $f$ \ is
 continuous (see, e.g., Jacod and Shiryaev \cite[VI.1.6]{JSh}).
In a similar way,
 \[
   \left\| f^{(n)}\left(\frac{j}{n}\right)
           - f^{(n)}\left(\frac{j-1}{n}\right) \right\|
   \leq \omega_T(f, n^{-1}) + 2 \sup_{t \in [0,T]} \|f^{(n)}(t) - f(t)\| 
 \]
 for all \ $j \in \{1, \ldots, n\}$.
\ By \eqref{rate},
 \[
   \sum_{j=1}^{\nt} \bigl\| \bm_\bxi^{(\nt-j)r} - \bPi \bigr\|
   \leq \sum_{j=1}^{\nT} c \kappa^{\nt-j}
   \leq \frac{c}{1 - \kappa} .
 \]
Further,
 \begin{align*}
  &\left\| \left(\bm_\bxi^{\nt r} -\bPi\right)
           \Bigg[ \sum_{\ell=i}^r \bm_{\bxi}^{\ell-i} f^{(n)}_\ell(0)
                  + \sum_{\ell=1}^{i-1}
                     \bm_{\bxi}^{\ell-i+r} f^{(n)}_\ell(0) \Bigg] \right\| \\
  &\leq \left\| \left(\bm_\bxi^{\nt r} -\bPi\right)
           \Bigg[ \sum_{\ell=i}^r \bm_{\bxi}^{\ell-i} (f^{(n)}_\ell(0) - f_\ell(0))
                  + \sum_{\ell=1}^{i-1}
                     \bm_{\bxi}^{\ell-i+r}
                     (f^{(n)}_\ell(0) - f_\ell(0)) \Bigg] \right\| \\
  &\quad
        + \left\| \left(\bm_\bxi^{\nt r} -\bPi\right)
                  \Bigg[ \sum_{\ell=i}^r \bm_{\bxi}^{\ell-i} f_\ell(0)
                         + \sum_{\ell=1}^{i-1}
                            \bm_{\bxi}^{\ell-i+r} f_\ell(0) \Bigg] \right\| \\
  &\leq c K \sum_{\ell=1}^r \|f^{(n)}_\ell(0) - f_\ell(0)\| ,
 \end{align*}
 since \ $f \in C$ \ implies
 \[
   \left(\bm_\bxi^{\nt r} - \bPi\right)
   \Bigg[ \sum_{\ell=i}^r \bm_{\bxi}^{\ell-i} f_\ell(0)
          + \sum_{\ell=1}^{i-1} \bm_{\bxi}^{\ell-i+r} f_\ell(0) \Bigg]
   = \bzero .
 \]
Indeed, by \ $f \in C$ \ and \eqref{Pimr}, we obtain
 \begin{align*}
  \bm_\bxi^{\nt r}
  \Bigg[ \sum_{\ell=i}^r \bm_{\bxi}^{\ell-i} f_\ell(0)
         + \sum_{\ell=1}^{i-1} \bm_{\bxi}^{\ell-i+r} f_\ell(0) \Bigg]
  &=\bm_\bxi^{\nt r} \bPi
    \Bigg[ \sum_{\ell=i}^r \bm_{\bxi}^{\ell-i} f_\ell(0)
           + \sum_{\ell=1}^{i-1} \bm_{\bxi}^{\ell-i+r} f_\ell(0) \Bigg] \\
  &=\bPi
    \Bigg[ \sum_{\ell=i}^r \bm_{\bxi}^{\ell-i} f_\ell(0)
           + \sum_{\ell=1}^{i-1} \bm_{\bxi}^{\ell-i+r} f_\ell(0) \Bigg] .
 \end{align*}
Using that \ $f^{(n)} \lu f$ \ as \ $n \to \infty$, \ we have
 \ $\Psi^{(n)}(f^{(n)}) \lu \psi(f)$ \ as \ $n \to \infty$.
\ Thus we conclude \ $C \subset C_{\Psi, \, (\Psi^{(n)})_{n \in \NN}}$.

By the definition of a strong solution (see, e.g., Jacod and Shiryaev
 \cite[Definition 2.24, Chapter III]{JSh}), \ $\bcM$ \ has almost sure
 continuous sample paths, so we have \ $\PP(\bcM \in C) = 1$.
\ Consequently, by Lemma \ref{Conv2Funct}, we obtain
 \ $\bcY^{(n)} = \Psi^{(n)}(\bcM^{(n)}) \distr \Psi(\bcM) \distre \bcX$ \ as
 \ $n \to \infty$.
\proofend

\section{Appendix}
\label{Appendix}

In the proof of Theorem \ref{main} we will use some facts about the first and
 second order moments of the sequences \ $(\bX_k)_{k \in \ZZ_+}$ \ and
 \ $(\bM_k)_{k \in \NN}$.

\begin{Lem}\label{Moments}
Let \ $(\bX_k)_{k\in\ZZ_+}$ \ be a $p$-type branching process with immigration.
Suppose that \ $\bX_0 = \bzero$, \ $\EE(\|\bxi_{1,1,i}\|^2) < \infty$ \ for all
 \ $i \in \{ 1, \dots, p \}$ \ and \ $\EE (\|\bvare_1\|^2) < \infty$.
\ Then
 \begin{gather}
  \EE(\bX_k)
  = \sum_{j=0}^{k-1} \bm_{\bxi}^j \bm_{\bvare}, \label{mean}\\
  \var(\bX_k)
  = \sum_{j=0}^{k-1} \bm_{\bxi}^j \bV_{\bvare} (\bm_{\bxi}^\top)^j
    + \sum_{j=0}^{k-2}
       \bm_{\bxi}^j
       \sum_{\ell=0}^{k-j-2}
        \left[ (\bm_{\bxi}^\ell \bm_{\bvare}) \odot \bV_\bxi \right]
        (\bm_{\bxi}^\top)^j . \label{var}     
 \end{gather}
If, in addition, \ $(\bX_k)_{k\in\ZZ_+}$ \ is a critical indecomposable $p$-type
 branching process with immigration and the offspring mean matrix \ $\bm_\bxi$
 \ has the form \eqref{partition_m}, then, for all \ $k \in \NN$, \ we have
 \ $\EE(\bM_k \mid \cF_{rk-r}) = \bzero$, \ $\EE(\bM_k) = \bzero$, \ and
 \begin{align} \label{Mcond} 
  \EE(\bM_k \bM_k^\top \mid \cF_{rk-r})
  &=\bigoplus_{\ell = 1}^r
     \left\{ \left[ \bm_{\bxi}^{r-\ell} \bX_{rk-r}
                    + \sum_{j=1}^{r-\ell} \bm_{\bxi}^{j-1} \bm_\bvare \right] 
             \odot \bV_{\bxi} + \bV_{\bvare} \right\} , \\
  \EE(\bM_k \bM_k^\top)
  &=\bigoplus_{\ell = 1}^r
     \left\{ \left[ \bm_{\bxi}^{r-\ell} \EE(\bX_{rk-r})
                    + \sum_{j=1}^{r-\ell} \bm_{\bxi}^{j-1} \bm_\bvare \right] 
             \odot \bV_{\bxi} + \bV_{\bvare} \right\} . \label{Cov}
 \end{align}
\end{Lem}

\noindent
\textbf{Proof.}
We have already proved \eqref{mean}, see \eqref{EXk}.
The equality
 \ $\bM_{k,\ell} = \bX_{rk-\ell+1} - \EE(\bX_{rk-\ell+1} \mid \cF_{rk-\ell})$
 \ clearly implies \ $\EE(\bM_{k,\ell} \mid \cF_{rk-\ell}) = \bzero$, \ thus
 \ $\EE(\bM_{k,\ell} \mid \cF_{rk-r})
    = \EE[\EE(\bM_{k,\ell} \mid \cF_{rk-\ell}) \mid \cF_{rk-r}]
    = \bzero$,
 \ and hence \ $\EE(\bM_k) = \bzero$.
\ The proof of \eqref{var} can be found in Isp\'any and Pap \cite{IP}.
By \eqref{MBPI(d)} and \eqref{Mk},
 \begin{equation}\label{Mdeco}
  \begin{aligned}
   \bM_{k,\ell}
   &= \bX_{rk-\ell+1} - \EE(\bX_{rk-\ell+1} \mid \cF_{rk-\ell})
    = \bX_{rk-\ell+1} - \bm_{\bxi} \bX_{rk-\ell} - \bm_{\bvare} \\
   &= \sum_{i=1}^p \sum_{j=1}^{X_{rk-\ell,i}}
       ( \bxi_{rk-\ell+1,j,i} - \EE(\bxi_{rk-\ell+1,j,i}) )
       + ( \bvare_{rk-\ell+1} - \EE(\bvare_{rk-\ell+1}) ) .
  \end{aligned}
 \end{equation}
For each \ $k \in \NN$ \ and \ $\ell \in \{1, \ldots, p\}$, \ the random
 vectors
 \[
   \big\{\bxi_{rk-\ell+1,j,i} - \EE(\bxi_{rk-\ell+1,j,i}) , \,
          \bvare_{rk-\ell+1} - \EE(\bvare_{rk-\ell+1})
          : j \in \NN , \, i \in \{1, \dots, p\}
   \big\}
 \]
 are independent of each others, independent of \ $\cF_{rk-\ell}$, \ and have
 zero mean, hence
 \[
   \EE(\bM_{k,\ell} \bM_{k,\ell}^\top \mid \cF_{rk-\ell})
   = \sum_{i=1}^p X_{rk-\ell,i} \bV_{\bxi_i} + \bV_\bvare
   = [\bX_{rk-\ell} \odot \bV_\bxi] + \bV_\bvare .
 \]
By the tower rule and by \eqref{mart},
 \[
   \EE(\bM_{k,\ell} \bM_{k,\ell}^\top \mid \cF_{rk-r})
   = [\EE(\bX_{rk-\ell} \mid \cF_{rk-r}) \odot \bV_\bxi] + \bV_\bvare
   = \left[ \bm_{\bxi}^{r-\ell} \bX_{rk-r}
            + \sum_{j=1}^{r-\ell} \bm_{\bxi}^{j-1} \bm_\bvare \right] 
     \odot \bV_{\bxi}
     + \bV_{\bvare} .
 \]
If \ $j, \ell \in \{1, \ldots, p\}$ \ with \ $j < \ell$ \ then, again by the
 tower rule,
 \[
   \EE(\bM_{k,j} \bM_{k,\ell}^\top \mid \cF_{rk-r})
   = \EE(\bM_{k,j} \EE(\bM_{k,\ell} \mid \cF_{rk-\ell})^\top \mid \cF_{rk-r})
   = \bzero
 \]
 since \ $\EE(\bM_{k,\ell} \mid \cF_{rk-\ell}) = \bzero$, \ and similarly, if
 \ $j, \ell \in \{1, \ldots, p\}$ \ with \ $j > \ell$ \ then
 \ $\EE(\bM_{k,j} \bM_{k,\ell}^\top \mid \cF_{rk-r}) = \bzero$, \ thus we conclude
 \eqref{Mcond}, and hence, \eqref{Cov}.
\proofend

\begin{Lem}\label{EEX}
Let \ $(\bX_k)_{k\in\ZZ_+}$ \ be a critical indecomposable $p$-type branching
 process with immigration.
Suppose that \ $\bX_0 = \bzero$, \ $\EE(\|\bxi_{1,1,i}\|^4) < \infty$ \ for all
 \ $i \in \{ 1, \dots, p \}$ \ and \ $\EE (\|\bvare_1\|^4) < \infty$.
\ Then
 \begin{gather*}
  \EE(\|\bX_k\|) = \OO(k) , \qquad
  \EE(\|\bX_k\|^2) = \OO(k^2) , \qquad
  \EE(\|\bM_k\|) = \OO(k^{1/2}), \qquad
  \EE(\|\bM_k\|^4) = \OO(k^2) .
 \end{gather*}
\end{Lem}

\noindent
\textbf{Proof.}
By \eqref{mean},
 \[
   \|\EE(\bX_k)\|
   \leq \sum_{j=0}^{k-1}
           \|\bm_{\bxi}^j\| \cdot \|\bm_{\bvare}\|
   = \OO(k) ,
 \]
 since
 \begin{equation}\label{C}
  C_\bxi := \sup_{j \in \ZZ_+} \| \bm_{\bxi}^j \| < \infty .
 \end{equation}
Indeed, write \ $j \in \ZZ_+$ \ in the form \ $j = rk + i$ \ with
 \ $k \in \ZZ_+$ \ and \ $i \in \{0, \ldots, r-1\}$.
\ Then
 \ $\| \bm_{\bxi}^j \| = \| \bm_{\bxi}^{rk + i} \|
    \leq \| \bm_{\bxi}^{rk} \| \| \bm_{\bxi}^i \|
    \leq (c + \|\bPi\|) \max_{i \in \{0, \ldots, r-1\}} \| \bm_{\bxi}^i \|
    =: C_\bxi < \infty$,
 \ since \eqref{rate} implies
 \ $\| \bm_{\bxi}^{rk} \| \leq \| \bm_{\bxi}^{rk} - \bPi \| + \|\bPi\|
    \leq c + \|\bPi\|$. 

We have
 \[
   \EE(\|\bX_k\|^2) = \EE\bigl[\tr(\bX_k \bX_k^\top)\bigr]
   = \tr(\var(\bX_k))
     + \tr\bigl[ \EE(\bX_k) \EE(\bX_k)^\top \bigr] ,
 \]
 where
 \ $\tr\bigl[ \EE(\bX_k) \EE(\bX_k)^\top \bigr]
    = \|\EE(\bX_k)\|^2 \leq \bigl[\EE(\|\bX_k\|)\bigr]^2
    = \OO(k^2)$.
\ Moreover, \ $\tr(\var(\bX_k)) = \OO(k^2)$.
\ Indeed, by \eqref{var} and \eqref{C},
 \begin{align*} 
  \|\var(\bX_k)\|
  &\leq \|\bV_{\bvare}\|
        \sum_{j=0}^{k-1} \| \bm_{\bxi}^j \|^2
        + \|\bm_{\bvare}\| \cdot \|\bV_\bxi\|
          \sum_{j=0}^{k-2}
           \|\bm_{\bxi}^j\|^2 \sum_{\ell=0}^{k-j-2} \|\bm_{\bxi}^\ell\| \\
  &\leq C_\bxi^2 \|\bV_{\bvare}\| k
        + C_\bxi^3 \|\bm_{\bvare}\| \cdot \|\bV_\bxi\| k^2 ,
 \end{align*}
 where \ $\|\bV_\bxi\| := \sum_{i=1}^p \|\bV_{\bxi_i}\|$, \ hence we obtain
 \ $\EE(\|\bX_k\|^2) = O(k^2)$.

We have
 \begin{align*}
  \EE(\|\bM_k\|)
  &\leq \sqrt{\EE(\|\bM_k\|^2)}
   = \sqrt{\EE\bigl[ \tr (\bM_k \bM_k^\top) \bigr]} \\
  &= \sqrt{\tr\left[ \bigoplus_{i = 1}^r
                      \left\{ \left[ \bm_{\bxi}^{r-i} \EE(\bX_{rk-r})
                                     + \sum_{j=1}^{r-i}
                                        \bm_{\bxi}^{j-1} \bm_\bvare \right] 
                              \odot \bV_{\bxi} + \bV_{\bvare} \right\} \right]} ,
 \end{align*}
 hence we obtain \ $\EE(\|\bM_k\|) = \OO(k^{1/2})$ \ from
 \ $\EE(\|\bX_k\|) = \OO(k)$. 

By \eqref{Mdeco}, 
 \[
   \|\bM_{k,\ell}\|
   \leq \|\bvare_{rk-\ell+1} - \EE(\bvare_{rk-\ell+1})\|
        + \sum_{i=1}^p
           \Biggl\| \sum_{j=1}^{X_{rk-\ell,i}}
                    ( \bxi_{rk-\ell+1,j,i} - \EE(\bxi_{rk-\ell+1,j,i}) ) \Biggr\| ,
 \]
 hence
 \begin{align*}
  \EE(\|\bM_{k,\ell}\|^4)
  &\leq (p+1)^3 \EE(\|\bvare_1 - \EE(\bvare_1)\|^4) \\
  &\quad
        + (p+1)^3 
          \sum_{i=1}^p
           \EE\Biggl(\Biggl\| \sum_{j=1}^{X_{rk-\ell,i}}
                               ( \bxi_{rk-\ell+1,j,i} - \EE(\bxi_{rk-\ell+1,j,i}) )
                     \Biggr\|^4
              \Biggr) .
 \end{align*}
Here
 \begin{multline*}
  \EE\Biggl(\Biggl\| \sum_{j=1}^{X_{rk-\ell,i}}
                      ( \bxi_{rk-\ell+1,j,i} - \EE(\bxi_{rk-\ell+1,j,i}) ) \Biggr\|^4
            \Biggr) \\
   \begin{aligned}
    &= \EE\Biggl[\Biggl(\sum_{m=1}^p
                         \Biggl( \sum_{j=1}^{X_{rk-\ell,i}}
                                ( \xi_{rk-\ell+1,j,i,m} - \EE(\xi_{rk-\ell+1,j,i,m}) )
                       \Biggr)^2
                \Biggr)^2\Biggr] \\
    &\leq p
          \sum_{m=1}^p
           \EE\Biggl[\Biggl(\sum_{j=1}^{X_{rk-\ell,i}}
                            ( \xi_{rk-\ell+1,j,i,m} - \EE(\xi_{rk-\ell+1,j,i,m}) )
                     \Biggr)^4
              \Biggr] ,
  \end{aligned}
 \end{multline*}
 where
 \begin{multline*}
  \EE\Biggl[\Biggl(\sum_{j=1}^{X_{rk-\ell,i}}
                   ( \xi_{rk-\ell+1,j,i,m} - \EE(\xi_{rk-\ell+1,j,i,m}) ) \Biggr)^4
            \, \Bigg| \, \cF_{kr-\ell}\Biggr] \\
  = X_{rk-\ell,i} \EE[( \xi_{1,1,i,m} - \EE(\xi_{1,1,i,m}) )^4]
     + X_{rk-\ell,i} (X_{rk-\ell,i} - 1)
       \bigl(\EE[( \xi_{1,1,i,m} - \EE(\xi_{1,1,i,m}) )^2]\bigr)^2 
 \end{multline*}
 with
 \ $\bigl(\EE[( \xi_{1,1,i,m} - \EE(\xi_{1,1,i,m}) )^2]\bigr)^2
    \leq \EE[( \xi_{1,1,i,m} - \EE(\xi_{1,1,i,m}) )^4]$,
 \ hence, by \ $X_{rk-\ell,i} \geq 0$ \ and
 \ $X_{rk-\ell,i} (X_{rk-\ell,i} - 1) \geq 0$,
 \[
   \EE\Biggl[\Biggl(\sum_{j=1}^{X_{rk-\ell,i}}
                    ( \xi_{rk-\ell+1,j,i,m} - \EE(\xi_{rk-\ell+1,j,i,m}) ) \Biggr)^4
      \Biggr]
   \leq \EE[( \xi_{1,1,i,m} - \EE(\xi_{1,1,i,m}) )^4]
        \EE[(X_{rk-\ell,i})^2] .
 \]
Consequently, \ $\EE(\|\bX_k\|^2) = \OO(k^2)$ \ implies
 \ $\EE(\|\bM_k\|^4) = \OO(k^2)$.
\proofend

Next we recall a result about convergence of random step processes towards a
 diffusion process, see Isp\'any and Pap \cite[Corollary 2.2]{IspPap}.

\begin{Thm}\label{Conv2DiffThm}
Let \ $\bgamma : \RR_+ \times \RR^p \to \RR^{p \times q}$ \ be a continuous
 function.
Assume that uniqueness in the sense of probability law holds for the SDE
 \begin{equation}\label{SDE}
  \dd \, \bcU_t
  = \gamma (t, \bcU_t) \, \dd \bcW_t ,
  \qquad t \in \RR_+,
 \end{equation}
 with initial value \ $\bcU_0 = \bu_0$ \ for all \ $\bu_0 \in \RR^p$, \ where
 \ $(\bcW_t)_{t \in \RR_+}$ \ is a $q$-dimensional standard Wiener process.
Let \ $(\bcU_t)_{t \in \RR_+}$ \ be a solution of \eqref{SDE} with initial value
 \ $\bcU_0 = \bzero$.

For each \ $n \in \NN$, \ let \ $(\bU^{(n)}_k)_{k\in\NN}$ \ be a sequence of
 $p$-dimensional martingale differences with respect to a filtration
 \ $(\cF_k)_{k \in \ZZ_+}$.
\ Let
 \[
   \bcU^{(n)}_t := \sum_{k=1}^{\nt} \bU^{(n)}_k \, ,
   \qquad t \in \RR_+, \quad n \in \NN .
 \]
Suppose \ $\EE \big( \|\bU^{(n)}_k\|^2 \big) < \infty$ \ for all
 \ $n, k \in \NN$.
\ Suppose that for each \ $T > 0$,
 \begin{enumerate}
  \item[\textup{(i)}]
   $\sup\limits_{t\in[0,T]}
     \left\| \sum\limits_{k=1}^{\nt}
              \EE\Bigl[ \bU^{(n)}_k (\bU^{(n)}_k)^\top \mid \cF_{k-1} \Bigr]
             - \int_0^t
                \bgamma(s,\bcU^{(n)}_s) \bgamma(s,\bcU^{(n)}_s)^\top
                \dd s \right\|
         \stoch 0$,\\
  \item[\textup{(ii)}]
   $\sum\limits_{k=1}^{\lfloor nT \rfloor}
     \EE \big( \|\bU^{(n)}_k\|^2 \bbone_{\{\|\bU^{(n)}_k\| > \theta\}}
               \bmid \cF_{k-1} \big)
    \stoch 0$
   \ for all \ $\theta > 0$,
 \end{enumerate}
 where \ $\stoch$ \ denotes convergence in probability.
Then \ $\bcU^{(n)} \distr \bcU$ \ as \ $n \to \infty$.
\end{Thm}

Now we recall a version of the continuous mapping theorem.

For functions \ $f$ \ and \ $f_n$, \ $n \in \NN$, \ in \ $\DD(\RR_+, \RR^p)$,
 \ we write \ $f_n \lu f$ \ if \ $(f_n)_{n \in \NN}$ \ converges to \ $f$
 \ locally uniformly, i.e., if \ $\sup_{t \in [0,T]} \|f_n(t) - f(t)\| \to 0$ \ as
 \ $n \to \infty$ \ for all \ $T > 0$.
\ For measurable mappings \ $\Phi : \DD(\RR_+, \RR^p) \to \DD(\RR_+, \RR^q)$
 \ and \ $\Phi_n : \DD(\RR_+, \RR^p) \to \DD(\RR_+, \RR^q)$, \ $n \in \NN$, \ we
 will denote by \ $C_{\Phi, (\Phi_n)_{n \in \NN}}$ \ the set of all functions
 \ $f \in \CC(\RR_+, \RR^p)$ \ for which \ $\Phi_n(f_n) \to \Phi(f)$ \ whenever
 \ $f_n \lu f$ \ with \ $f_n \in \DD(\RR_+, \RR^p)$, \ $n \in \NN$.

\begin{Lem}\label{Conv2Funct}
Let \ $(\bcU_t)_{t \in \RR_+}$ \ and \ $(\bcU^{(n)}_t)_{t \in \RR_+}$, \ $n \in \NN$,
 \ be \ $\RR^p$-valued stochastic processes with c\`adl\`ag paths such that
 \ $\bcU^{(n)} \distr \bcU$.
\ Let \ $\Phi : \DD(\RR_+, \RR^p) \to \DD(\RR_+, \RR^q)$ \ and
 \ $\Phi_n : \DD(\RR_+, \RR^p) \to \DD(\RR_+, \RR^q)$, \ $n \in \NN$, \ be
 measurable mappings such that there exists \ $C \subset C_{\Phi,(\Phi_n)_{n\in\NN}}$
 \ with \ $C \in \cD_\infty(\RR_+, \RR^p)$ \ and \ $\PP(\bcU \in C) = 1$.
\ Then \ $\Phi_n(\bcU^{(n)}) \distr \Phi(\bcU)$.
\end{Lem}

Lemma \ref{Conv2Funct} can be considered as a consequence of Theorem 3.27 in
 Kallenberg \cite{K}, and we note that a proof of this lemma can also be found
 in Isp\'any and Pap \cite[Lemma 3.1]{IspPap}.

\begin{Lem}\label{measurability}
The mappings \ $\Psi^{(n)}$, \ $n \in \NN$, \ and \ $\Psi$ \ defined in Section
 \ref{Proof_Conv_bX} are measurable.
\end{Lem}

\noindent
\textbf{Proof.}
Continuity of \ $\Psi$ \ follows from the characterization of convergence in
 \ $\DD(\RR_+, (\RR^p)^r)$, \ see, e.g., Ethier and Kurtz
 \cite[Proposition 3.5.3]{EK}, thus we obtain measurability of \ $\Psi$.

In order to prove measurability of \ $\Psi^{(n)}$, \ first we localize it.
For each \ $N \in \NN$, \ consider the stopped mapping
 \ $\Psi^{(n,N)} : \DD(\RR_+, (\RR^p)^r) \to \DD(\RR_+, (\RR^p)^r)$ \ given by
 \ $\Psi^{(n,N)}(f)(t) := \Psi^{(n)}(f)(t \land N)$ \ for
 \ $f \in \DD(\RR_+, (\RR^p)^r)$, \ $t \in \RR_+$, \ $n, N \in \NN$.
\ Obviously, \ $\Psi^{(n,N)}(f) \to \Psi^{(n)}(f)$ \ in \ $\DD(\RR_+, (\RR^p)^r)$
 \ as \ $N \to \infty$ \ for all \ $f \in \DD(\RR_+, (\RR^p)^r)$, \ since for
 all \ $T > 0$ \ and \ $N \geq T$ \ we have
 \ $\Psi^{(n,N)}(f)(t) = \Psi^{(n)}(f)(t)$, \ $t \in [0,T]$, \ and hence
 \ $\sup_{t \in [0,T]} \|\Psi^{(n,N)}(f)(t) - \Psi^{(n)}(f)(t)\| \to 0$ \ as
 \ $N \to \infty$.
\ Consequently, it suffices to show measurability of \ $\Psi^{(n,N)}$ \ for all
 \ $n, N \in \NN$.
\ We can write \ $\Psi^{(n,N)} = \Psi^{(n,N,2)} \circ \Psi^{(n,N,1)}$, \ where the
 mappings \ $\Psi^{(n,N,1)} : \DD(\RR_+, (\RR^p)^r) \to ((\RR^p)^r)^{nN+1}$ \ and
 \ $\Psi^{(n,N,2)} : ((\RR^p)^r)^{nN+1} \to \DD(\RR_+, (\RR^p)^r)$ \ are defined by
 \begin{align*}
   \Psi^{(n,N,1)}(f)
   &:= \left( f(0), f\left( \frac{1}{n} \right), f\left( \frac{2}{n} \right),
              \dots, f(N) \right) , \\
   \Psi^{(n,N,2)}(\bx_0, \bx_1, \dots, \bx_{nN})(t)
   &:= \bm_{\bxi}^{\nt r}
       \Bigg[ \sum_{\ell=i}^r \bm_{\bxi}^{\ell-i} \bx_{0,\ell}
             + \sum_{\ell=1}^{i-1} \bm_{\bxi}^{\ell-i+r} \bx_{0,\ell} \Bigg] \\
  &\quad\:
      + \sum_{j=1}^\nt
         \bm_{\bxi}^{(\nt-j)r} 
         \Bigg[ \sum_{\ell=i}^r
                  \bm_{\bxi}^{\ell-i}
                  \left( \bx_{j,\ell} - \bx_{j-1,\ell}
                         + \frac{1}{nr} \bm_{\bvare} \right) \\
  &\phantom{\quad\: + \sum_{j=1}^\nt \bm_{\bxi}^{(\nt-j)r} \Bigg[}
      + \sum_{\ell=1}^{i-1}
         \bm_{\bxi}^{\ell-i+r}
         \left( \bx_{j-1,\ell} - \bx_{j-2,\ell}
                + \frac{1}{nr} \bm_{\bvare} \right) \Bigg]
 \end{align*}
 for \ $f \in \DD(\RR_+, (\RR^p)^r)$, \ $t \in \RR_+$,
 \ $\bx = (\bx_0, \bx_1, \dots, \bx_{nN}) \in ((\RR^p)^r)^{nN+1}$,
 \ $n, N \in \NN$.
\ Measurability of \ $\Psi^{(n,N,1)}$ \ follows from Ethier and Kurtz
 \cite[Proposition 3.7.1]{EK}.
Next we show continuity of \ $\Psi^{(n,N,2)}_n$ \ by checking
 \ $\sup_{t \in [0,T]}
     \|\Psi^{(n,N,2)}(\bx^{(k)})(t) - \Psi^{(n,N,2)}(\bx)(t)\| \to 0$
 \ as \ $k \to \infty$ \ for all \ $T > 0$ \ whenever \ $\bx^{(k)} \to \bx$
 \ in \ $((\RR^p)^r)^{nN+1}$.
 \ This convergence follows from the estimates
 \begin{align*}
  &\sup_{t \in [0,T]} \|\Psi^{(n,N,2)}(\bx^{(k)})(t) - \Psi^{(n,N,2)}_n(x)(t)\| \\
  &\leq C_\bxi^2 \sum_{\ell=1}^r \|\bx^{(k)}_{0,\ell} - \bx_{0,\ell}\|
        + C_\bxi^2
          \sum_{j=1}^{\lfloor n (T \land N) \rfloor}
           \Biggl[ \sum_{\ell=i}^r
                    \bigl( \|\bx^{(k)}_{j,\ell} - \bx_{j,\ell}\|
                           + \|\bx^{(k)}_{j-1,\ell} - \bx_{j-1,\ell}\| \bigr) \\
  &\phantom{\leq C_\bxi^2 \sum_{\ell=1}^r \|\bx^{(k)}_{0,\ell} - \bx_{0,\ell}\|
        + C_\bxi^2
          \sum_{j=1}^{\lfloor n (T \land N) \rfloor}
           \Biggl[}
                   + \sum_{\ell=1}^{i-1}
                      \bigl( \|\bx^{(k)}_{j-1,\ell} - \bx_{j-1,\ell}\|
                             + \|\bx^{(k)}_{j-2,\ell} - \bx_{j-2,\ell}\| \bigr)
           \Biggr] .
 \end{align*}
We obtain measurability of both \ $\Psi^{(n,N,1)}$ \ and \ $\Psi^{(n,N,2)}$,
 \ hence we conclude measurability of \ $\Psi^{(n,N)}$.
\proofend

\begin{Lem}\label{closed}
The subset \ $\CC(\RR_+, (\RR^p)^r) \subset \DD(\RR_+, (\RR^p)^r)$ \ is closed,
 thus measurable, i.e.,
 \ $\CC(\RR_+, (\RR^p)^r) \in \cD_\infty(\RR_+, (\RR^p)^r)$.
\end{Lem}

\noindent
\textbf{Proof.}
The complement \ $\DD(\RR_+, (\RR^p)^r) \setminus \CC(\RR_+, (\RR^p)^r)$ \ is
 open.
Indeed, each function
 \ $f \in \DD(\RR_+, (\RR^p)^r) \setminus \CC(\RR_+, (\RR^p)^r)$ \ is
 discontinuous at some point \ $t_f \in \RR_+$, \ and, by the definition of the
 metric of \ $\DD(\RR_+, (\RR^p)^r)$, \ there exists \ $r_f > 0$ \ such that
 all \ $g \in \DD(\RR_+, (\RR^p)^r)$ \ is discontinuous at the point
 \ $t_f \in \RR_+$ \ whenever the distance of \ $g$ \ from \ $f$ \ is less than
 \ $r_f$.
\ Consequently, the set
 \ $\DD(\RR_+, (\RR^p)^r) \setminus \CC(\RR_+, (\RR^p)^r)$ \ is the union of
 open balls with center
 \ $f \in \DD(\RR_+, (\RR^p)^r) \setminus \CC(\RR_+, (\RR^p)^r)$ \ and radius
 \ $r_f$.
\proofend

\end{document}